\documentclass[openany]{book}
\usepackage{Wiley-AuthoringTemplate}
\usepackage[authoryear]{natbib}
\usepackage{amsmath}
\usepackage{subcaption}
\usepackage{algorithm}
\usepackage{algpseudocode}
\usepackage{amsmath,bm}
\usepackage{caption}
\usepackage{tabularx}
\usepackage{cleveref}
\usepackage{array}
\usepackage{amsmath,bm}
\usepackage{siunitx}
\usepackage{amsmath}
\usepackage{epstopdf}
 \usepackage{graphicx}
 \usepackage{booktabs}
 \usepackage{rotating}
\usepackage{tabularx}
\usepackage{hyperref}
 \usepackage{tikz}
 \usepackage{comment}

 \usepackage{setspace}
 \usepackage{longtable}
 \usepackage{xtab}
\usepackage{mathtools,nccmath}
 \usepackage{nicefrac}
 \usepackage{float}
 \usepackage{footnote}
 \usepackage{xcolor}
 \usepackage{titlesec}
\renewcommand{\thesection}{\arabic{section}}
\renewcommand{\thesubsection}{\thesection.\arabic{subsection}}

\usepackage{setspace}  
  \newcommand{\set}{\mathcal}
\begin{document}

\title{\textbf{\Huge 
    \vspace{10pt}A Tutorial on Multi-time Scale \\ 
    Optimization Models and \vspace{10pt}
    \\ Algorithms}
}

\author[1]{Asha Ramanujam}
\author*[1]{Can Li}
\address[1]{\orgdiv{Davidson School of Chemical Engineering},
\orgname{Purdue University},
\postcode{47907},  \countrypart{Indiana},
\city{West Lafayette}, \street{480 W.Stadium Ave}, \country{USA}}
\address*{Corresponding Author: Can Li; \email{canli@purdue.edu}}
\maketitle
\begin{abstract}{Abstract}
Systems across different industries consist of interrelated processes and decisions in different time scales including long-time decisions and short-term decisions. To optimize such systems, the most effective approach is to formulate and solve multi-time scale optimization models that integrate various decision layers. In this tutorial, we provide an overview of multi-time scale optimization models and review the algorithms used to solve them. We also discuss the metric \textit{Value of the Multi-scale Model} (VMM) introduced to quantify the benefits of using multi-time scale optimization models as opposed to sequentially solving optimization models from high-level to low-level. Finally, we present an illustrative example of a multi-time scale capacity expansion planning model and showcase how it can be solved using some of the algorithms (\url{https://github.com/li-group/MultiScaleOpt-Tutorial.git}). This tutorial serves as both an introductory guide for beginners with no prior experience and a high-level overview of current algorithms for solving multi-time scale optimization models, catering to experts in process systems engineering.
\end{abstract}
\begin{objectives}
\item Explain multi-time scale optimization models and their structure
\item Discuss \textit{Value of the Multi-scale Model} (VMM) as metric to assess the advantages of using multi-time scale optimization models
\item Present an overview of the algorithms used in literature to solve multi-time scale optimization models 
\item Examine an example of a multi-time scale capacity expansion planning model and solve it using some of the algorithms discussed
\end{objectives}

\section{Introduction}
Large process systems can consist of different interrelated processes. These systems involve decisions taken at different time scales including long-time, high-level planning decisions (e.g., supply chain design, capacity planning, and production targets) on a yearly, monthly, or weekly basis and short-time, low-level scheduling and control decisions on a daily or hourly scale or even in seconds \citep{Shin2019Multi-timescaleProgramming,Maravelias2009IntegrationOpportunities,Allen2023ASystems,Subramanian2013IntegrationManagement}.

Optimizing a system requires determining both high-level and low-level decisions. One approach is to sequentially solve models that involve high-level and low-level decisions. In this approach, a high-level problem is first formulated and solved to derive the high-level decisions. These high-level decisions are then fixed in the low-level problems to obtain the low-level decisions. This process creates a one-way, top-down communication within the hierarchical structure.

However, making high-level decisions without any input or feedback from low-level problems often leads to infeasible or suboptimal solutions. This occurs because the decision layers are interconnected and influence one another. Therefore, to ensure a feasible and optimal solution, it is crucial to integrate feedback from low-level problems into the high-level decision-making process. This integration acknowledges the interdependencies between decision layers and enhances overall system performance.

This integration can be achieved by formulating and solving a multi-time scale optimization model that integrates decision processes taking place at different time scales \citep{Shin2019Multi-timescaleProgramming,Maravelias2009IntegrationOpportunities,Allen2023ASystems,Subramanian2013IntegrationManagement}.  Multi-time scale optimization models have been formulated to solve problems in various domains, such as plant production and scheduling \citep{Biondi2017OptimizationApproach}, electricity market participation \citep{Dowling2017AParticipation} and environmental management\citep{ Tabrizi2018IntegratedImpacts}. In particular, an area where multi-time scale optimization models are significant is in the decarbonization and electrification of chemical process systems \citep{Ramanujam2023DistributedMicrogrid,Kim2024Multi-periodUncertainty}. Electrification entails powering chemical processes with electricity which can be sourced from renewable green energy rather than fossil fuels. However, these renewable resources have substantial spatial and temporal variations, and electricity prices also fluctuate over time. To fully exploit the economic possibilities of the entire system, it is crucial to account for these fluctuations in the design of electrified chemical plants and the renewable energy sources that fuel them. The optimal approach includes integrating the planning and scheduling of these plants and resources into a multi-timescale optimization model \citep{Ramanujam2023DistributedMicrogrid}. 

This tutorial explores multi-time scale optimization models and the algorithms designed to solve them. It is structured as follows: in section \ref{MSOM}, we give a representation for a general multi-time scale optimization model. We describe \textit{Value of the Multi-scale Model} (VMM), a metric used to quantify the benefits of using a multi-time scale optimization model in section \ref{sec:VMM}. In section \ref{sec:algorithms}, we discuss the algorithms used in literature to solve multi-time scale optimization models. We provide an example of a multi-time scale optimization model in section \ref{sec:example} and discuss the implementation of a few of the algorithms on the example. Finally, in section \ref{sec:conclusion}, we draw the conclusions.

\section{Multi-Time Scale Optimization Models} \label{MSOM}

As discussed before, multi-time scale optimization models take decisions at different time scales. We divide the decision variables into two sets: high-level decision variables and low-level decision variables. The high-level decision variables denoted by $\bm{x}$ pertain to the long time scales and low-level decision variables denoted by $\bm{y}_s$ encompass all the other shorter time scales defined for subperiods $s$.  
Multi-time scale optimization models can be represented in the following way:
\begin{subequations}
    \label{FP}
    \begin{align}
      \text{MM} =  \min_{\bm{x},\bm{y}_s
       } & \quad f(\bm{x}) + \sum_{s \in S}w_s q(\bm{x}, \bm{y}_s,\bm{\theta}_s) & \\
        \text{s.t.} \quad & g(\bm{x}) \leq 0 \label{FP1} & \\
        & h(\bm{x}, \bm{y}_s, \bm{\theta}_s) \leq 0 \quad \forall s \in S \label{FP2} & \\
        & \Tilde{h}(\bm{x}, \bm{y}_1,\bm{y}_2,...,\bm{y}_{|S|}, \bm{\theta}) \label{FP3} \leq 0 & 
    \end{align}
\end{subequations}
The notation used for the above formulation is given in Table \ref{Notation:mm}.

To understand the formulation, consider an integrated design and scheduling problem. The problem can be formulated as a multi-time scale optimization model. In the multi-time scale model, the design decisions are the high-level decisions \(\bm{x}\) and the scheduling decisions are the low-level decisions \(\bm{y}_s\) for various subperiods \(s\). The model can have various constraints which can be classified in the following way: 
\begin{itemize}
    \item Equation \eqref{FP1} represents constraints solely on the high-level decisions. Examples include constraints on the bounds of the design decisions, such as the maximum possible capacity for reactor vessels.
    \item Constraints concerning the high-level decisions and the low-level decisions of a single subperiod \(s\) are represented by equation \eqref{FP2}. Examples include constraints on the bounds of scheduling decisions like maximum possible production.
     \item Constraints concerning the low-level decisions for multiple single subperiods and possibly the high-level decisions are represented by equation \eqref{FP3}. Such constraints are often called complicating constraints.  Examples include constraints that combine the production of multiple subperiods to satisfy the demand for the entire duration.
    
\end{itemize}

This method of writing the constraints is general and will help us further to understand whether different algorithms can be applied to different multi-time scale optimization model in Section \ref{sec:algorithms}. 

\begin{table}
\caption{Notation used in Equation \eqref{FP}\label{Notation:mm}}{%
\begin{tabular}{l l}
\toprule
Notation & Meaning \\
\midrule
$\bm{x}$ & high-level decisions \\
$s$ & subperiod \\
$S$ & set of subperiods \\
$w_s$ & probability/weight of subperiod $s$ \\
$\bm{\theta}_s$ &  parameters associated with subperiod $s$ \\
$\bm{y}_s$ &  low-level decision variables associated with subperiod $s$ \\
$\bm{\theta}$ &  complete set of parameters \\
$f$ & cost function associated with the high-level decisions \\
$q$ & cost function associated with the low-level decisions \\
\(g\) & constraints pertaining only to the high-level decisions \\ 
$h$ & constraints concerning only a single subperiod \(s\) \\
$\Tilde{h}$ & constraints connecting the decisions in different subperiods \\
\botrule
\end{tabular}}{}
\end{table}

\textbf{Differences with 2-stage stochastic programming}:
2-stage stochastic programming  models have a formulation similar to \eqref{FP} with \(\bm{x}\) representing the first-stage variables, and \(\bm{y}_s\) representing the second-stage variables for scenario \(s.\)  The major difference between these models and multi-time scale models is that these models do not have complicating constraints that connect decisions in different subperiods i.e., no \( \Tilde{h}.\) However, multi-time scale models can have complicating constraints.

\section{Value of the Multi-scale Model (VMM)} \label{sec:VMM}
\begin{figure}[H]
\centering
 \includegraphics[width = 5cm]{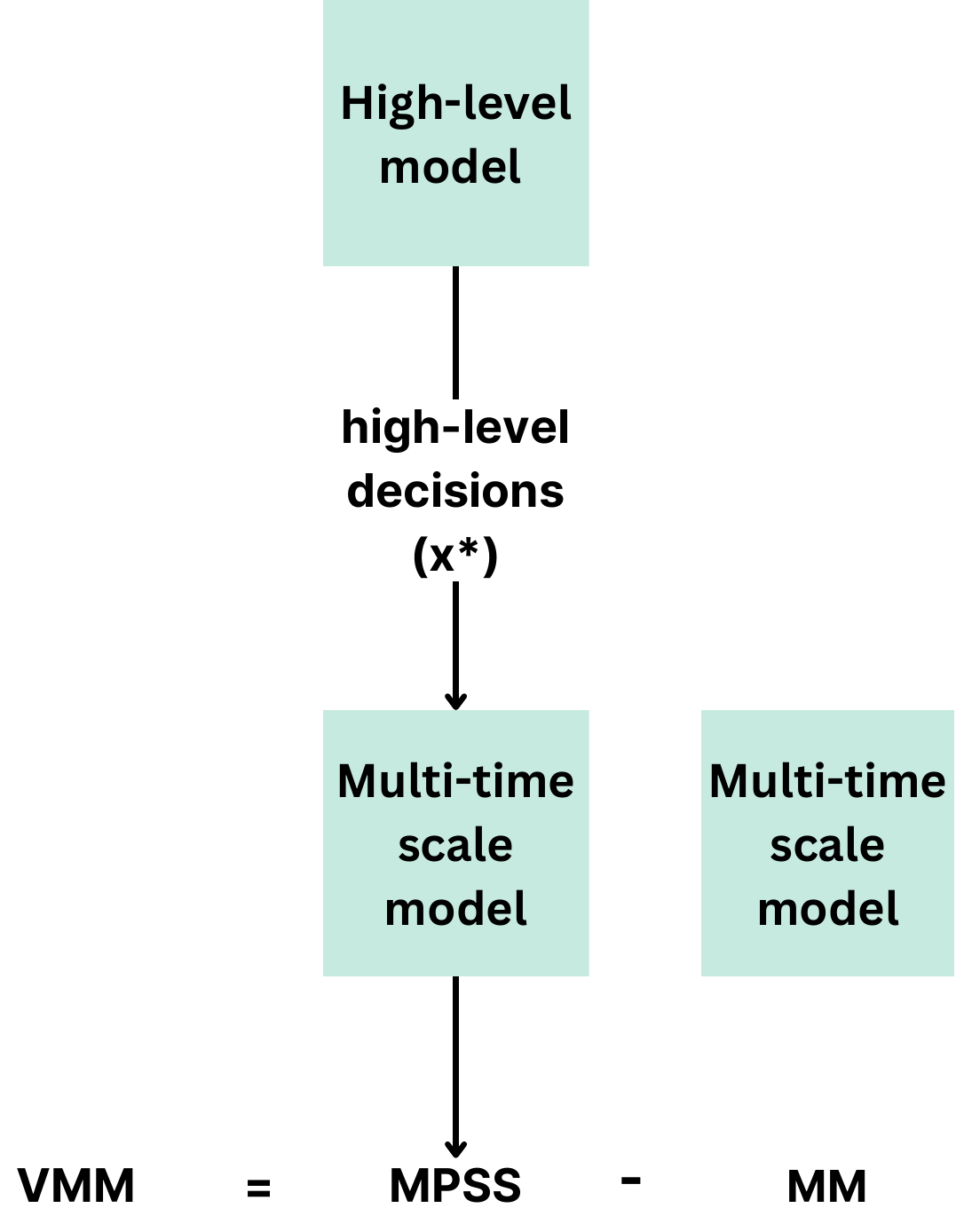}
\caption{\textit{Value of the Multi-scale Model}}
\label{fig:VMM}
\end{figure}

While solving a multi-time scale algorithm might be very difficult as compared to a traditional one-way top-down approach, it provides a more economic solution. To assess the benefits of using a multi-time scale optimization model instead of the traditional one-way top-down approach, we can use the metric \textit{Value of the Multi-scale Model} (VMM) defined by \cite{Ramanujam2023DistributedMicrogrid}. The VMM is the difference between the objective derived from the one-way top-down approach (MPSS) and the optimal objective from solving the multi-time scale model (MM) as shown in Figure \ref{fig:VMM}. The VMM helps evaluate the enhancements in system efficiency and feasibility gained by integrating decision-making across various time scales.

We now define the VMM more formally. Let's start by introducing some notations. The optimal objective  obtained from solving the multi-time scale optimization model will be referred to as the performance of the multi-time scale model (MM).

Next, consider a simpler, high-level model that is used in the traditional top-down approach. This model can be formulated as follows:

\begin{subequations}
\label{UP}
\begin{align}
\text{SM} = \min_{\bm{x},\bm{z}} \quad & f(\bm{x}) + Q(\bm{x},\bm{z},\bm{\Tilde{\theta}}) & \\
\text{s.t.} \quad & g(\bm{x}) \leq 0 & \\
& H(\bm{x}, \bm{z}, \bm{\Tilde{\theta}}) \leq 0 &
\end{align}
\end{subequations}

In this formulation, $\bm{z}$ serves as a surrogate for the low-level decisions, and $\bm{\Tilde{\theta}}$ as a surrogate parameter for the original parameters $\bm{\theta}_s$. Similarly, \(Q\) is a surrogate for the function \(q\) while \(H\) encompasses surrogates for both \(h\) and \(\Tilde{h}\) i.e., the equivalent function and constraints in terms of \(\bm{z}\) and $\bm{\Tilde{\theta}}$.

These high-level decisions are then fixed in the multi-time scale model to obtain the low-level decisions. The model shown below is the multi-time scale model with these decisions fixed:
\begin{subequations}
\label{LP}
\begin{align}
\text{MPSS} = \min_{\bm{y}s} \quad & f(\bm{x}^*) + \sum_{s \in S}w_s q(\bm{x}^*, \bm{y}_s,\bm{\theta}_s) & \\
\text{s.t.} \quad & h(\bm{x}^*, \bm{y}_s, \bm{\theta}_s) \leq 0 \quad \forall s \in S & \\
& \Tilde{h}(\bm{x}^*, \bm{y}_1,\bm{y}_2,...,\bm{y}_{|S|}, \bm{\theta}) \leq 0 &
\end{align}
\end{subequations}

where $\bm{x}^*$ signifies the high-level decisions obtained from solving the high-level model in \eqref{UP}. The model \eqref{LP} can also be referred to as the low-level model.

The Multi-scale Performance of the Single-scale Solution (MPSS) is defined as the optimal objective obtained from solving \eqref{LP} as shown in Figure \ref{fig:VMM}. In case, the \eqref{LP} is infeasible, we set MPSS to be a very large positive value like \(10^{10}.\)

Finally, the \textit{Value of the Multi-scale Model} (VMM) is calculated as the difference between the Multi-scale Performance of the Single-scale Solution (MPSS) and the performance of the multi-time scale optimization model (MM).

 \[VMM = MPSS-MM\]

This metric quantifies the benefits of using a multi-time scale model over the traditional one-way top-down approach. A higher VMM indicates greater improvements in system efficiency and feasibility, justifying the additional complexity of the multi-time scale optimization.

The VMM metric has been proposed by taking inspiration from the Value of the Stochastic Solution (VSS) introduced by the stochastic programming community to quantify the benefits of stochastic solutions as compared to optimal solutions from deterministic models \citep{Birge2011IntroductionProgramming}.
Computing VSS involves solving an \textit{Expected Value} problem for one nominal scenario corresponding to the expected value of $\bm{\theta}_{s}$, i.e., $\bm{\bar{\theta}} = \sum_{s\in \set{S}}w_{s} \bm{\theta}_{s}$ 

\begin{subequations} \label{eq:expected}
  \begin{align}
  \text{EV} = \min_{\bm{x}, \bm{\bar{y}}} f(\bm{x}) + q(\bm{x}, \bm{\bar{y}}, \bm{\bar{\theta}}), \\ \text{s.t.}\quad g(\bm{x})\leq 0 \\ \quad h(\bm{x}, \bm{\bar{y}}, \bm{\bar{\theta}}) \leq 0
\end{align}  
\end{subequations}

In the above formulation, \(\bm{\bar{y}}\) represents the low-level decisions for the nominal scenario. Typically, the size of the expected value model \eqref{eq:expected} is larger than the size of the high-level model used to compute VMM \eqref{UP} since it still involves multiple time scales. 

 The performance of EV is then evaluated by substituting the solution $\bm{x^{EV}}$ obtained from \eqref{eq:expected} in the original full-space problem with the optimal objective being called the \textit{expected result of using the EV solution} (EEV) \citep{Birge2011IntroductionProgramming}.
\begin{subequations}
\label{LP_EEV}
\begin{align}
\text{EEV} = \min_{\bm{y}s} \quad & f(\bm{x^{EV}}) + \sum_{s \in S}w_s q(\bm{x^{EV}}, \bm{y}_s,\bm{\theta}_s) & \\
\text{s.t.} \quad & h(\bm{x^{EV}}, \bm{y}_s, \bm{\theta}_s) \leq 0 \quad \forall s \in S & 
\end{align}
\end{subequations}

 VSS is then defined by
\[
    VSS = EEV - MM.
\]

The difference between VMM and VSS is the way in which the high-level decisions are obtained for the top-down approach. While VMM uses an aggregated model with a lesser number of time scales involved, VSS uses a nominal variable involving the same number of time scales. Thus VMM shows the value of having variables in multiple time scales, while VSS illustrates the value of having multiple scenarios/subperiods.  

 \section{Algorithms to solve multi-time scale optimization models} \label{sec:algorithms}
 Multi-time scale optimization models can be of different sizes (from thousands of variables and constraints to billions of variables and constraints) and different structures. Different algorithms have been proposed to solve multi-time scale optimization models. We classify these algorithms in the following way:
 \begin{enumerate}
     \item \hyperref[fpmethods]{Full-space methods}
     \item \hyperref[decompmethods]{Decomposition algorithms}
     \item \hyperref[metaheuristic]{Metaheuristic algorithms}  
     \item \hyperref[mathheuristic]{Matheuristic algorithms} 
     \item \hyperref[datadriven]{Data-driven methods}
     \item \hyperref[PAMSO]{Parametric Autotuning Multi-time Scale Optimization (PAMSO) algorithm }
 \end{enumerate}

 In the next few subsections, we will be discussing these algorithms. Note that in this chapter, when we refer to solving a model with an algorithm, we refer to the process of generating good feasible solutions for the problem. Apart from full-space methods and conventional decomposition algorithms, the other algorithms used to address multi-time scale optimization models generally do not provide theoretical guarantees of optimality or valid bounds (on their own) for assessing solution quality.
 \subsection{Full-space methods} \label{fpmethods}
 One method to solve multi-time scale optimization models is to solve them directly using exact optimization techniques. In other words, we can solve the proposed MILP/MINLP integrated models directly using solvers like Gurobi, CPLEX, or BARON. For example, \cite{Kopanos2009Multi-SiteIndustries} proposed an MILP to solve production planning and scheduling/batching for multi-site batch process industries. A case study involving 20 final products in 2 different production plants was solved using the CPLEX 11.0 solver. The MILP for the case study involved 1,444 continuous variables, 5,256 binary variables, and  6,744 constraints and was solved in 430 seconds.

 While this can done for small-size problems, it becomes very difficult and time-consuming to use these solvers to solve problems with millions of variables and constraints. This is because MILP and MINLP problems are NP-hard problems in general making the problem more intractable as the problem sizes increase.

\subsection{Decompositon Algorithms}\label{decompmethods}
 Decomposition algorithms involve breaking the multi-time scale optimization model into smaller subproblems. They can be used on problems moderately large (up to tens of millions of variables \Citep{Li2022Mixed-integerSystems}) and tackle specific types of problems.  The most commonly used decomposition algorithms are the following:

 \begin{enumerate}
     \item \hyperref[bilevel]{Bi-level decomposition} 
     \item \hyperref[Benders]{Benders decomposition} 
     \item \hyperref[para:lag]{Lagrangian decomposition} 
     \item \hyperref[bp]{Branch and price}  
 \end{enumerate}

 We discuss the decomposition algorithms below. Furthermore, we discuss in detail the implementation of some of the algorithms to solve a simplified linear multi-time scale optimization model \eqref{Decomp:FP}. We use a simplified problem as it helps illustrate some of the algorithm's classical versions, which are then suitably modified to solve more complex multi-time scale optimization models. The model of the simplified problem is shown below: 
 \begin{subequations}
 \label{Decomp:FP}
      \begin{align}
\min \quad & c^T \bm{x}+\sum_{s \in \set{S}}q_{s}^T \bm{y}_{s}&  \label{Decomp:FP1}\\
\text { s. t.} \quad & A \bm{x}=b, \quad \bm{x} \geq 0 &  \label{Decomp:FP2}\\
& T_{s} \bm{x}+W_{s} \bm{y}_{s}=h_{s}, \quad \bm{y}_{s} \geq 0,  \quad \forall s\in \set{S}&  \label{Decomp:FP3}
\end{align}
 \end{subequations}

In the above representation, \(\bm{x}\) represents the high-level decisions and \(\bm{y}_s\) represents the low-level decisions for subperiod \(s.\) Furthermore, \eqref{Decomp:FP2} signifies the constraints only on the high-level decisions and \eqref{Decomp:FP3} denotes the constraints connecting the high-level and low-level decisions. In the model we consider, we do not have complicating constraints i.e., no \( \Tilde{h}.\)  This is because some of the decomposition algorithms do not work on problems with such constraints.

As discussed before, without \(\Tilde{h},\) the model \eqref{Decomp:FP} can be treated as a two-stage stochastic linear program, with the subperiods being the scenarios, \(\bm{x}\) being the first stage variables and \(\bm{y}_s\) being the second stage variables. In this work, we treat this problem as a multi-time scale optimization problem as opposed to a stochastic program.

\subsubsection{Bi-level decomposition} \label{bilevel}
\begin{figure}[H]
\centering
\includegraphics[width=10cm]{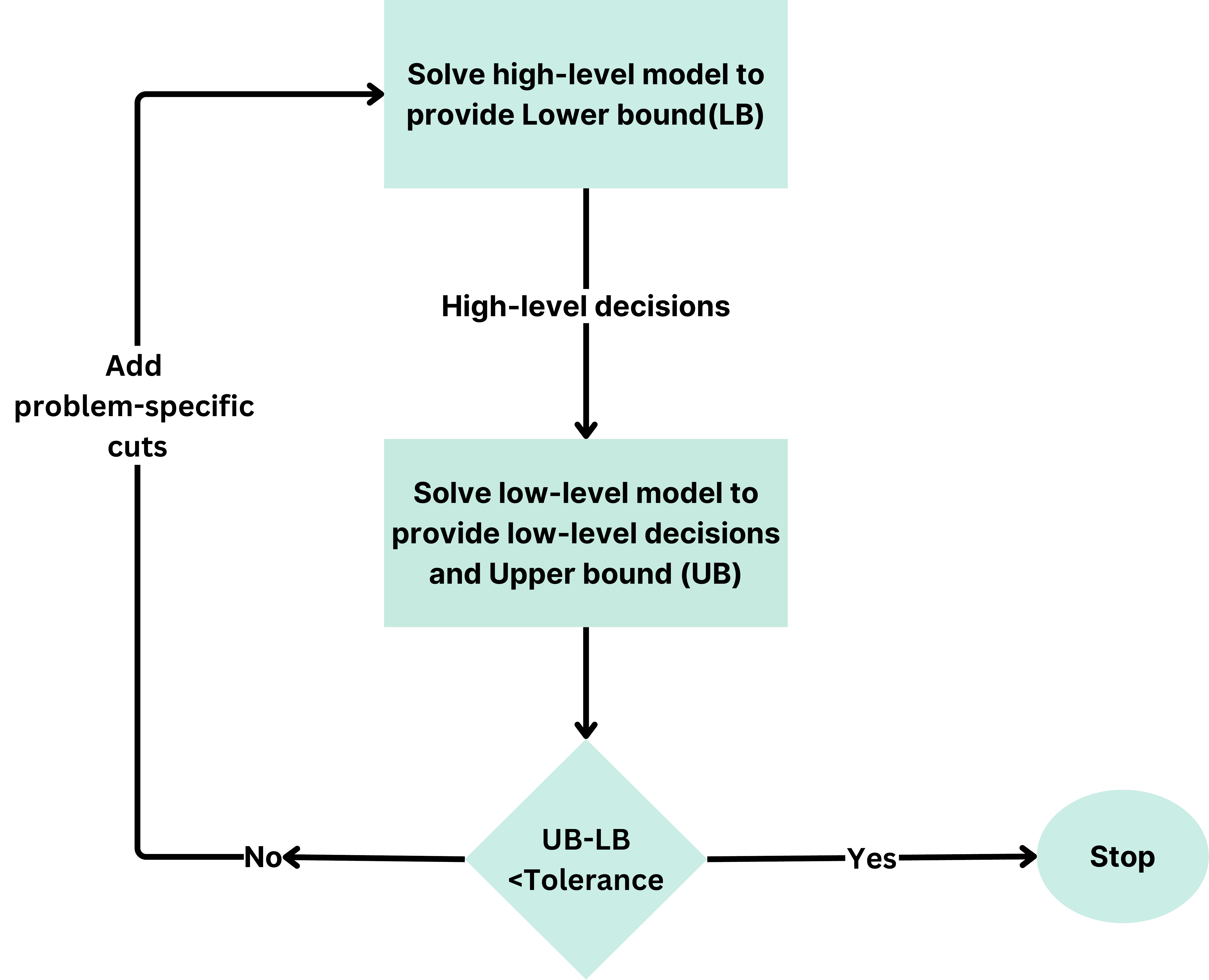}
\caption{Schematic for Bi-level decomposition}
\label{fig:bilevel}
\end{figure}

Bi-level decomposition involves decomposing the original model into a high-level model and a low-level model with the performance of the high-level model improved using feedback from the low-level model. The high-level model, typically a relaxation of the optimization model provides bounds for the model and some of the high-level decisions. The low-level model provides a feasible solution taking into account the high-level decisions obtained from the high-level model as shown in Figure \ref{fig:bilevel}.  
The algorithm involves iteratively solving the high-level and low-level models to reduce the difference in the objective values obtained from the models. Problem-specific logic and integer cuts based on the solution from the low-level problem are added to the high-level model to improve the bounds and high-level decisions generated. 

The steps involved are shown in algorithm \ref{alg:bi-level}.

\begin{algorithm}[H]
\caption{Bi-level decomposition algorithm}
\label{alg:bi-level}
\begin{algorithmic}[1]
\State \textbf{Initialize:} Set the Upper Bound (UB) $\leftarrow \infty$ and Lower Bound (LB) $\leftarrow -\infty$.
\State Decompose the problem into a high-level model and a low-level model.
\While{True}  
    \State Solve the high-level model to obtain high-level decisions and update the Lower Bound (LB).
    \State Fix the high-level decisions in the low-level model and solve the low-level model to obtain low-level decisions and update the Upper Bound (UB).
    \If{$\text{UB} - \text{LB} \leq \text{Tolerance}$}
        \State \textbf{Stop:} Convergence achieved.
    \Else
        \State Add problem-specific cuts to the high-level model.
    \EndIf
\EndWhile
\end{algorithmic}
\end{algorithm}

Bi-level decomposition algorithms have been used in literature to solve multi-time scale optimization models. For example, \cite{Erdirik-Dogan2008SimultaneousLines} presented a multi-period MILP integrating the planning and scheduling of single-stage multi-product continuous plants with parallel units.  To solve the model with longer time horizons, a bi-level decomposition algorithm was implemented decomposing the problem into a high-level and low-level model.  The high-level model formulated was based on a relaxation of the original model and provided an upper bound for the profit, the assignments of products to available equipment during each time period and the sequencing of the products. The low-level model excluded the products that were not selected by the high-level problem for each unit at each time period in the full-space model and provided detailed scheduling. Integer cuts were passed as feedback from the low-level model to the high-level model to exclude supersets of previously obtained feasible configurations. The algorithm was tested on several examples and generally outperformed full-space methods using commercial solvers by achieving faster solutions. The algorithm required only a single iteration to solve these examples. This demonstrates that the high-level model was well-formulated.

The bi-level decomposition algorithm can be extended to the multi-level decomposition algorithm that involves more than two levels. For instance,  \cite{Munawar2005IntegrationManufacturing}  presented a multi-level decomposition-based framework to solve the integration of planning and scheduling in a multi-site, multi-product plant. The framework was demonstrated for integrating planning and scheduling in the paper manufacturing industry. For small-sized problems, a 2-level decomposition algorithm was proposed with the high-level model deciding the simultaneous order allocation and sequencing at multiple sites and the low-level model minimizing the trim loss and sequencing of patterns. For medium to large-size problems, a 4-level decomposition algorithm was proposed with levels for order allocation, detailed scheduling, trim loss minimization, and sequencing of patterns.
\subsubsection{Dual-based decomposition algorithms}
In this part, we discuss the decomposition algorithms that involve the use of dual variables namely: Benders decomposition, Lagrangian decomposition, and Branch and price algorithm.

\begin{figure}[!ht]
\centering
\begin{subfigure}[b]{0.4\textwidth}
   \includegraphics[width = 5cm]{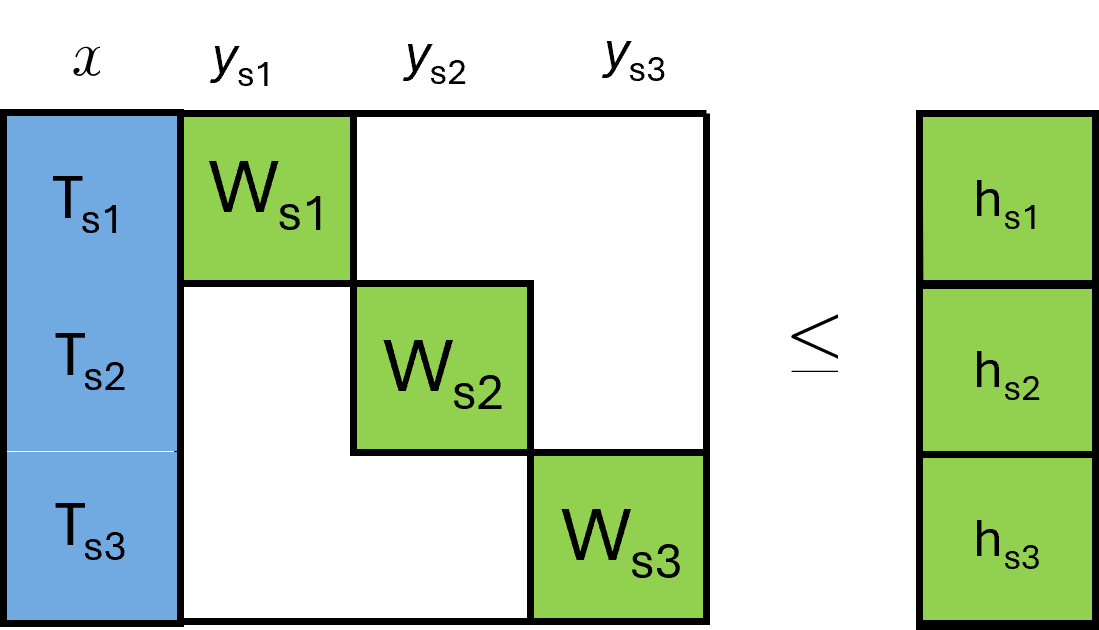}
   \caption{Benders decomposition}
    \label{fig:block:Benders}
\end{subfigure}
\hfill
\begin{subfigure}[b]{0.4\textwidth}
   \includegraphics[width = 5cm]{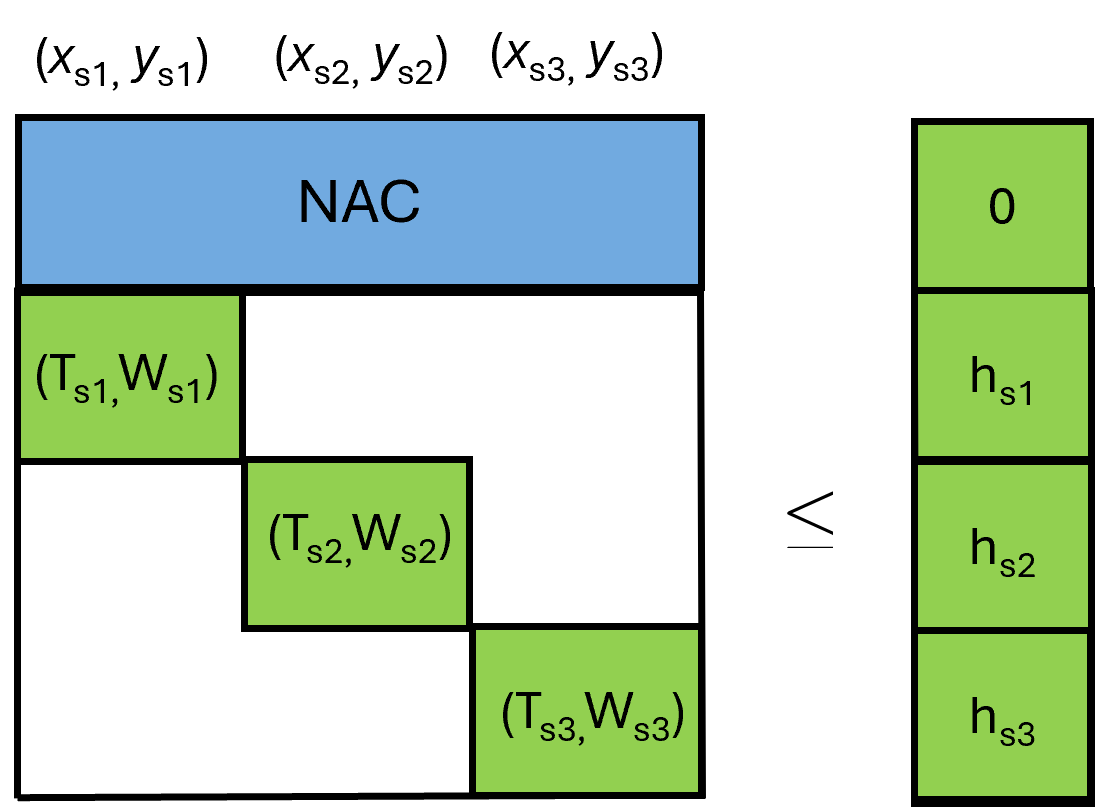}
    \caption{Lagrangian Decomposition and Dantzig Wolfe Decomposition}
    \label{fig:block:ld_cg}
\end{subfigure}
\caption{Block representation for decomposition algorithms}
\end{figure}

These decomposition methods solve problems with different structures. Benders decomposition solves problems with complicating variables like problems that have a representation as shown in Figure \ref{fig:block:Benders}. Complicating variables are those which, when fixed, make the remaining problem easier to solve in blocks. On the other hand, Dantzig Wolfe Decomposition (which serves as a key component of Branch and price algorithm) and  Lagrangian Decomposition solve problems with complicating constraints like problems that have a representation as shown in Figure \ref{fig:block:ld_cg}. Complicating constraints are constraints that when relaxed make the rest of the problem decomposable by blocks.  A specific type of complicating constraint is the Nonanticipativity Constraint (NAC), which ensures that variables in different blocks or subperiods are equalized.

\begin{figure}[H]
\centering
\includegraphics[width=12cm]{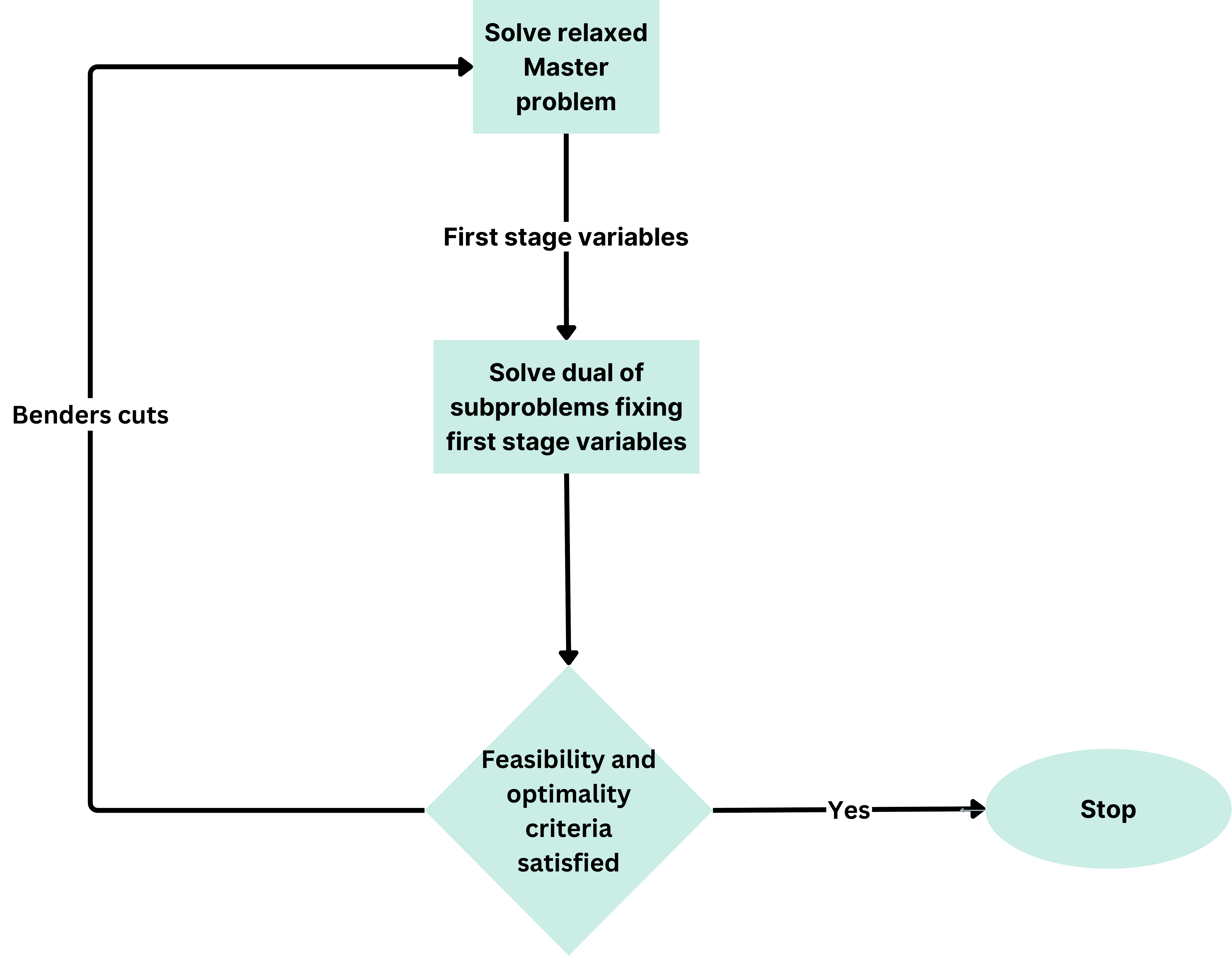}
\caption{Schematic for Classical Benders decomposition}
\label{fig:benders}
\end{figure}

We describe the algorithms in detail below.
\paragraph{Benders decomposition} \label{Benders}
Benders decomposition is a decomposition algorithm that can be used to solve multi-time scale optimization models with complicating variables. Benders decomposition involves splitting the problem into 2 sets of problems: a high-level master problem and low-level subproblems. The key idea is to iteratively solve these two problems, using information in the form of cuts derived from the subproblems to refine the master problem. This feedback loop continues until an optimal solution is found as shown in Figure \ref{fig:benders}. 

We now explain the classical version of Benders decomposition which requires that the subproblems be linear \Citep{Chu2013IntegratedApproach}. Consider solving the problem 
\eqref{Decomp:FP} using the Benders Decomposition.

 We first decompose the problem into a master problem  \eqref{MP_benders} and a set of subproblems \eqref{SP_benders} in the following way:

    \begin{subequations}
    \label{MP_benders}
    \begin{align}
      \text{MP} : & \min_{\bm{x},z}  c^T \bm{x}+\sum_{s \in \set{S}}z_s & \\
        \text{s.t.} \quad & A\bm{x} = b  & \\
        & \bm{x} \geq 0
    \end{align}
\end{subequations}

\begin{subequations}
    \label{SP_benders}
    \begin{align}
      \text{SP}_s :  & \min_{y} q^T_s \bm{y} & \\
        \text{s.t.} \quad & T_s\bm{x^*} + W_s\bm{y}  = h_s & \\
        & \bm{y} \geq 0 & 
    \end{align}
\end{subequations}

where \(\bm{x^*}\) represents the high-level decisions obtained from solving the master problem \eqref{MP_benders}.
If the solutions obtained from the subproblems are feasible and optimal we can stop. If that is not the case, we can improve the performance of the master problem by generating Benders cuts and updating the master problem.
To understand the optimality and feasiblity of the solution, we investigate the solution obtained from the subproblems. We need the subproblem to be feasible and finite. We can approach this by looking into the dual of the subproblem:

   \begin{subequations}
    \label{SP_bendersdual}
    \begin{align}
      \text{SP}^{dual}_s :  \max_{\bm{p}
       } \bm{p}^T(h_s-T_s\bm{x^*})  & \\
        \text{s.t.} \quad 
        \bm{p}^TW_s \leq q_s^T & 
    \end{align}
\end{subequations}
with \(p\) being the dual variables.

We assume the dual is feasible. This is always true when the primal is bounded. 

Three cases can occur for the optimization model corresponding to the dual of the subproblem \(s\):
\begin{itemize}
    \item \textbf{Case 1}: When the dual subproblem \eqref{SP_bendersdual} is unbounded.
    
    Here, the corresponding primal problem is infeasible. In such a case, there exists an extreme ray $\bm{w}_{s}^*$ in the polyhedron $\{\bm{p}|\bm{p}^T W_s \leq q_s^{T}\}$  such that $(\bm{w}_{s}^*)^{T}\left(h_s-T_s \bm{x^*}\right)>0.$ To handle this case, we generate a cut $(\bm{w}_{s}^*)^{T}\left(h_s-T_s x\right)\leq 0$ which we call the feasibility cut to the master problem \eqref{MP_benders}.
    \item \textbf{Case 2}: When the subproblem is feasible and bounded but the solution \(\bm{p}_s^*\)  is such that $z_s < \left(\bm{p}_s^*\right)^T\left(h_s-T_s x\right)$ where $\bm{p}_s^*$ is an extreme point of the polyhedron $X_s = \{\bm{p}|\bm{p}^T W_s \leq q_s^{T}\}$.
    
   The optimal dual lower bounds the optimal primal solution. Therefore, we need to have
    $$z_s\geq \max _{\bm{p} \in X_s}\left(\bm{p}\right)^T\left(h_s-T_s x\right)$$ If this is not ensured, we generate the cut $z_s\geq \left(\bm{p}_s^*\right)^T\left(h_s-T_s x\right)$ which we call the optimality cut to the master problem. 
    \item \textbf{Case 3}: When the dual subproblems are feasible and bounded and  $z_s\geq \left(\bm{p}_s^*\right)^T\left(h_s-T_s x\right)$ with \(\bm{p}_s^*\) being the solution of the dual subproblem.
    In such a case, we do not need to add any cut to the master problem.
\end{itemize}

When all the subproblems are feasible,bounded and provide solutions such that $z_s\geq \left(\bm{p}_s^*\right)^T\left(h_s-T_s x\right)$ with \(\bm{p}_s^*\) being the solution of the dual subproblem for subperiod \(s\), the feasibility and optimality criteria are satisfied. That is,  the solution derived from the master problem is optimal and we can terminate the algorithm.

The steps to solve a general linear problem with complicating variables are shown in Algorithm \ref{alg:benders}.
\begin{algorithm}[H]
\caption{Benders Decomposition Algorithm}
\label{alg:benders}
\begin{algorithmic}[1]
\State \textbf{Step 1:} Decompose the problem into a master problem \eqref{MP_benders} and a set of subproblems \eqref{SP_benders}.
\While{true}
    \State Solve the master problem \eqref{MP_benders} to obtain high-level decisions $\bm{x^*}$. \label{step2:benders}
    \State Using the high-level decisions obtained in Line \ref{step2:benders}, solve the dual of the subproblems \eqref{SP_benders}.
    \If{The feasibility and optimality criteria re satisfied}
        \State \textbf{Stop:} The current solution is optimal.
    \Else
        \State Generate Benders cuts and update the master problem \eqref{MP_benders}.
    \EndIf
\EndWhile
\end{algorithmic}
\end{algorithm}

While the master problem gives lower bounds due to its relaxed nature, summing up the subproblems gives upper bounds.

 There have been different extensions of Benders decomposition like logic-based Benders decomposition \Citep{Hooker2003Logic-basedDecomposition} where the subproblems do not have to be linear. These extensions can be used to solve problems like \eqref{FP} with complicating variables.  
 For example, \cite{Barzanji2020DecompositionProblem} formulated a mixed integer linear programming model for the integrated process planning and scheduling (IPSS) problem. They developed decomposition approaches based on the logic-based benders decomposition (LBBD) to solve the problem. The LBBD algorithm involved a master problem determining the process plan and operation-machine assignment, and a sub-problem that optimizes the sequencing and scheduling decisions. To improve convergence and enhance algorithm performance, the authors introduced relaxations for the optimal makespan objective function, a Benders optimality cut based on the critical path, and a faster heuristic for solving the sub-problem. The proposed algorithm was tested on 16 standard benchmark instances from the literature, consistently yielding either the optimal solution or improving the best-known solutions in all cases. This demonstrated the algorithm's superiority over existing state-of-the-art methods.

 Another example is the work by \cite{Chu2013IntegratedApproach}  which addressed the the integration of scheduling and dynamic optimization for batch chemical processes.  An MINLP dynamic optimization problem was formulated for the integration and the generalized Benders decomposition was adapted to solve the model. The master problem was a scheduling problem with variable processing times and processing costs, while the subproblem comprised a set of separable dynamic optimization problems for the processing units. The proposed decomposition algorithm was tested on case studies, with results demonstrating a reduction in computational effort by orders of magnitude compared to general-purpose MINLP solvers.

\begin{figure}[H]
\centering
\includegraphics[width=10cm]{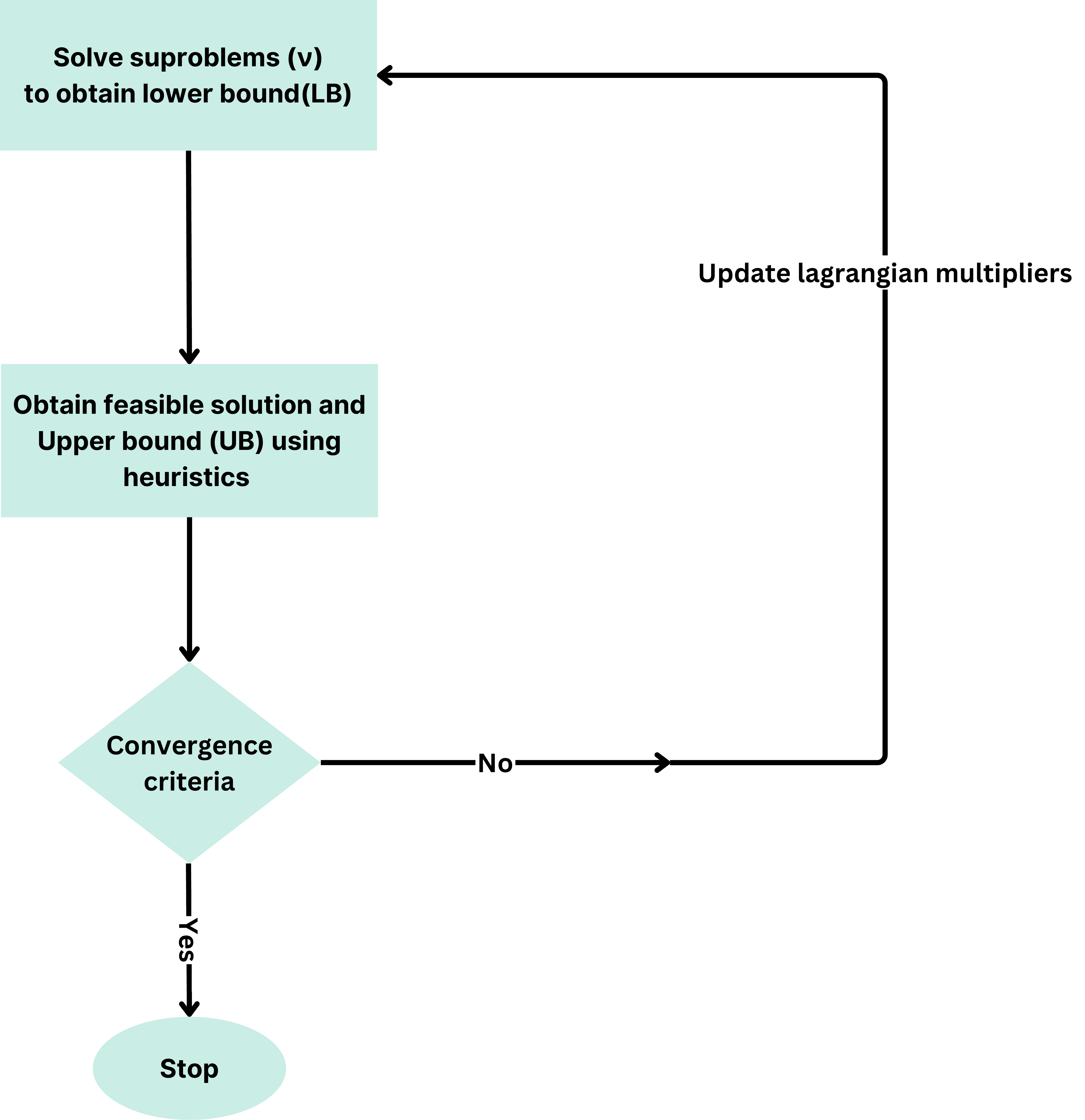}
\caption{Schematic for Lagrangian decomposition}
\label{fig:lag}
\end{figure}

 \paragraph{Lagrangian decomposition} \label{para:lag}
 Lagrangian decomposition \citep{Guignard1987LagrangeanBounds} involves relaxing specific constraints and penalizing their violation. This relaxation simplifies the original problem, allowing it to be decomposed into smaller subproblems, which can be solved independently The relaxation's performance improves through iterative updates based on the solutions of these subproblems as illustrated in Figure \ref{fig:lag}. 

Consider solving the optimization model \eqref{Decomp:FP} by decomposing it by subperiods. The structure of the problem originally looks like the figure \ref{fig:block:Benders}. 

One way to write this problem is by making copies the variable \(\bm{x}\) for each of the subperiod \(s\) i.e., we have a high-level decisions \(\bm{x}_s\) for each of the subperiod \(s.\) However, we need to ensure that the high-level decisions are same in all the subperiod \(s.\) Therefore, we add the constraints: \(\bm{x}_{1} = \bm{x}_{s} \quad \forall s \in \set{S},s \neq 1.\) These constraints are called Nonanticipativity constraints (NACs). With these modifications the original problem is written in the following way:
\begin{subequations} \label{lag1dec}
      \begin{align}
\min \quad & \sum_{s \in S}\left (\frac{c^T \bm{x}_s}{|\set{S}|}+q_{s}^T \bm{y}_{s} \right )&  \\
\text { s. t.} \quad & A \bm{x}_s=b,  \quad \bm{x}_s \geq 0 \quad  \forall s \in \set{S} & \\
& T_{s} \bm{x}_s+W_{s} \bm{y}_{s}=h_{s}, \quad \bm{y}_{s} \geq 0, \quad \forall s\in \set{S} &  \\
& \bm{x}_{1} = \bm{x}_{s} \quad \forall s \in \set{S},s \neq 1 &
\end{align}
 \end{subequations}

The structure of the problem now becomes like Figure \ref{fig:block:ld_cg}. It is still difficult to decompose the problem \eqref{lag1dec} due to the Nonanticipativity constraints. To address this, we relax these constraints and introduce a penalty for the deviation in high-level decisions, using weights known as Lagrangian multipliers. The process of relaxing and penalizing the violation of constraints is called ``dualizing the constraints'' \citep{Guignard1987LagrangeanBounds}. The relaxed problem is shown below:
\begin{subequations}\label{eq:lag_decomp_lagmul}
      \begin{align}
\min \quad & \sum_{s \in S}\left (\frac{c^T \bm{x}_s}{|\set{S}|}+q_{s}^T \bm{y}_{s} \right )+\sum_{s \in \set{S},s \neq 1} \bm{\nu}_s^T\left (\bm{x}_{1} - \bm{x}_{s}\right) &  \\
\text { s. t.} \quad & A \bm{x}_s=b,  \quad \bm{x}_s \geq 0 \quad  \forall s \in \set{S} & \\
& T_{s} \bm{x}_s+W_{s} \bm{y}_s=h_{s}, \quad  \bm{y}_s \geq 0, \forall s\in \set{S}& 
\end{align}
 \end{subequations}
 where \(\bm{\nu}_s\) are the Lagrangian multipliers.

We now decompose the model \eqref{eq:lag_decomp_lagmul} by subperiod s into the following subproblems:

\begin{subequations}
      \begin{align}
 D_1({\nu}):= \min_{\bm{x},\bm{y}} \quad & c^T \bm{x}+q_{1}^T \bm{y}+\sum_{s \in \set{S},s \neq 1} \bm{\nu}_s^T\bm{x} &  \\
\text { s. t.} \quad & A \bm{x}=b, \quad \bm{x} \geq 0  & \\
& T_{1} \bm{x}+W_{1} \bm{y}=h_{1} , \quad \bm{y} \geq 0 &  
\end{align}
 \end{subequations}

\begin{subequations}
      \begin{align}
 D_{s \in \set{S},s \neq 1}({\nu}):= \min_{\bm{x},\bm{y}} \quad & c^T \bm{x}+q_{s}^T \bm{y}-\bm{\nu}_s^T\bm{x} &  \\
\text { s. t.} \quad & A \bm{x}=b, \quad \bm{x} \geq 0  & \\
& T_{s} \bm{x}+W_{s} \bm{y}=h_{s}, \quad \bm{y} \geq 0 &
\end{align}
 \end{subequations}

In the above formulation,  \(\nu\) is the set of all the Lagrangian multipliers. These subproblems are relaxed as they are independent and do not explicitly take into account the equality of the high-level decisions between different subperiods. Thus, when the objectives of the subproblems are summed up they give a lower bound. We want to get the best lower bound by maximizing over all the possible Lagrangian multipliers. The best lower bound is  called the Lagrangian dual which is defined as \(v(LD) = \max_{\nu}  \sum_{s \in \set{S}}D_{s}({\nu}).\) Summing up all the objective functions of the subproblems, the Lagrangian dual can be written as the following way:

\[v(LD) = \max_{\nu}\min_{x_s^k,y_s^k}\sum_{s \in S}\left (\frac{c^T \bm{x}_s^k}{|\set{S}|}+q_{s}^T \bm{y}_s^k \right )+\sum_{s \in \set{S},s \neq 1} \bm{\nu}_s^T\left (\bm{x}_1^k - \bm{x}_s^k\right) \quad \forall (\bm{x}_s^k,\bm{y}_s^k) \in \set{K}_s, s \in S\]
 
 with \(X_s = \{\bm{x},\bm{y}|A \bm{x}=b, \bm{x} \geq 0,T_{s} \bm{x}+W_{s} \bm{y}=h_{s},\bm{y} \geq 0 \}, \text{ and }{K}_s = \text{Extreme points}( conv(X_s)).\) We assume that the subproblems are bounded. We use the extreme points of the convex hull because they define the boundaries of the feasible region.

The maximization-minimization problem expressed above can be written in the following way:

\begin{equation}
    \begin{aligned} \label{Lag_mp}
& \max_{\nu, \eta} \quad \eta \\
& \text { s.t. } \eta \leq \sum_{s \in S}\left (\frac{c^T \bm{x}_s^k}{|\set{S}|}+q_{s}^T \bm{y}_s^k \right )+\sum_{s \in \set{S},s \neq 1} \bm{\nu}_s^T\left (\bm{x}_1^k - \bm{x}_s^k\right) \quad \forall (\bm{x}_s^k,\bm{y}_s^k) \in \set{K}_s
\end{aligned}  
  \end{equation}

  The problem \eqref{Lag_mp}, can also be written in the following way:

  \begin{equation}
    \begin{aligned} \label{Lag_mp1}
& \max_{\nu, \eta} \quad \sum_{s \in S}\eta_s \\
& \text { s.t. } \eta_1 \leq \frac{c^T \bm{x}_1^k}{|\set{S}|}+q_{1}^T \bm{y}_1^k +\sum_{s \in \set{S},s \neq 1} \bm{\nu}_s^T\bm{x}_1^k \quad \forall (\bm{x}_1^k,\bm{y}_1^k) \in \set{K}_1\\
& \eta_s \leq\frac{c^T \bm{x}_s^k}{|\set{S}|}+q_{s}^T \bm{y}_s^k -\bm{\nu}_s^T \bm{x}_s^k \quad \forall (\bm{x}_s^k,\bm{y}_s^k) \in \set{K}_s,s \in \set{S},s \neq 1
\end{aligned}  
  \end{equation}

  \eqref{Lag_mp} is called the Lagrangian master problem with \eqref{Lag_mp1} being another way of writing the Lagrangian master problem. 

  The Lagrangian master problem \eqref{Lag_mp1} can be solved to get the optimal Lagrangian multipliers. It is very difficult to obtain all the extreme points and there can also be an exponential number of extreme points which makes the Lagrangian master problem intractable.  Therefore we can start with a subset of extreme points and iteratively update the master problem by obtaining more extreme points/ feasible solutions from the subproblems. That is, we can start with a set of multipliers, solve the subproblems using
the set multipliers, and update the master problem by adding the constraints 
$\eta_1 \leq \frac{c^T \bm{x}_1^k}{|\set{S}|}+q_{1}^T \bm{y}_1^k +\sum_{s \in \set{S},s \neq 1} \nu_s\bm{x}_1^k 
,\eta_s \leq\frac{c^T \bm{x}_s^k}{|\set{S}|}+q_{s}^T \bm{y}_s^k -\nu_s \bm{x}_s^k \quad \forall s \in \set{S},s \neq 1$ using the solution obtained \(\bm{x}_s^k,\bm{y}_s^k\) from the subproblems. We can then solve the updated master problem to obtain new multipliers. The constraint added is called the cutting plane and the method is called the cutting plane method \citep{Bertsekas1997NonlinearProgramming}.

Another way to find the Lagrangian multipliers is by using subgradients. In the above problem \eqref{Lag_mp}, \([\bm{x}_1-\bm{x}_s]_{s \in S,s \neq 1}\) represents a direction for the subgradient of the master problem. We can use the subgradients to solve the Lagrangian master problem

This method involves starting with an initial set of Lagrangian multipliers, solving the subproblems using the set multipliers, and updating the multipliers using the solutions obtained in the following way: 

\[\bm{\nu}_{k+1,s} = \bm{\nu}_{k,s}+\lambda_k(\bm{x}^k_1-\bm{x}^k_s) \quad \forall s \in \set{S},s \neq 1 \] 
where \(\bm{x}^k_1,\bm{x}^k_s\) are the solutions from solving the subproblems \( D_1, D_s\) respectively, \(\nu_{k,s}\) is the old Lagrangian multiplier, \(\nu_{k+1,s}\) is the new Lagrangian multiplier and \(\lambda_k\) is the step size. This method is called the subgradient method.

In practice a step size that has been proven to be effective is the following:

\[\lambda_{k} = \frac{v^*-v_k(LD)}{\sum_{s \in \set{S},s \neq 1}||\bm{x}^k_1-\bm{x}^k_s||^2}\]

where \(v_k(LD)\) is the Lagrangian dual obtained and \(v^*\) is the upper bound obtained from heuristics.

We can use the solutions obtained from the subproblems to obtain feasible solutions and correspondingly update upper bounds on the optimal objective. This iterative process is continued till the gap between the upper bound and lower bound is within a set a limit of tolerance or we obtain the same Lagrangian multipliers in 2 consecutive iterations.

The above can be extended to solve \eqref{FP} using the steps in Algorithm \ref{alg:lag}.

\begin{algorithm}[H]
\caption{Lagranigian Decomposition}
\label{alg:lag}
\begin{algorithmic}[1]
\State \textbf{Initialize:} Set the Upper Bound (UB) $\leftarrow \infty$, Lower Bound (LB) $\leftarrow -\infty$, k = 1.
\State Decompose the problem into independent subproblems by making copies of the high-level variables $\bm{x}$ and low-level decision variables $\bm{y}$ connected through the complicating constraints.
\State Start with an initial point for the Lagrangian multipliers $\nu_k$.
\While{true}
    \State Solve the subproblems using the Lagrangian multipliers $\nu_k$ and update the Lower Bound (LB).
    \State Use heuristics to obtain upper bounds on the optimal objective.
    \If{$\text{UB} - \text{LB} \leq \text{Tolerance}$ or  $k \geq 2,\bm{\nu}_k = \bm{\nu}_{k-1}$}
        \State \textbf{Stop:} Convergence achieved.
    \Else
        \State Update $k \rightarrow k+1$ and the Lagrangian multipliers $\bm{\nu}_k$. \label{algo:update}
        \State Go back to step 3.
    \EndIf
\EndWhile
\end{algorithmic}
\end{algorithm}

Line \ref{algo:update} of the algorithm can be done using the subgradient method or the cutting plane method.

Different extensions of Lagrangian decomposition/relaxation methods have been used to solve multi-time scale optimization models. For example, \cite{Shah2012IntegratedIndustry} provided a MILP formulation for an integrated planning and scheduling problem for the multisite batch facilities serving the global market. The augmented Lagrangian decomposition method was used to solve the formulation by dualizing the constraints linking planning and scheduling variables. Unlike classical Lagrangian decomposition, which typically uses a linear term in the Lagrangian function to enforce constraints, the augmented Lagrangian method incorporates a quadratic term. This quadratic term is used to penalize constraint violations, improving convergence properties and stability of the algorithm. To handle the non-separable cross-product terms in the augmented Lagrangian function, a diagonal approximation method was applied. The study of three examples by \cite{Shah2012IntegratedIndustry} demonstrated that the decomposition approach resulted in faster solution times compared to solving the full-space problem.

Progressive hedging \citep{Rockafellar1991ScenariosUncertainty} is another extension that relaxes the nonanticipativity constraints and solves the subproblems independently, with the algorithm penalizing deviations from the average values of decision variables. For instance, \cite{Peng2019AUncertainty} presented a novel progressive hedging-based algorithm to solve integrated planning and scheduling problems under demand uncertainty. A few strategies were proposed to accelerate the implementation of the algorithm as well as guarantee its convergence. The algorithm was tested on variants of a typical state–task network example as well as a real-world ethylene plant case, with results showing that it outperforms commercial solvers.

 \paragraph{Branch and price} \label{bp}
Column Generation methods are techniques used to solve large-scale linear programming (LP) problems that involve an overwhelming number of variables. Instead of solving the entire problem with all variables, these methods start with a restricted subset of variables, and iteratively improve this subset by solving subproblems. The subproblems identify new variables (columns) that can potentially improve the solution. This iterative approach continues until no further improvement can be made, effectively solving the full LP with a manageable number of variables.

The Dantzig-Wolfe Decomposition is a classical algorithm that employs column generation to solve linear programming problems with complicating constraints. Complicating constraints are those that connect variables across different subproblems, often creating dependencies that make solving the problem as a whole difficult. 

\begin{figure}[H]
\centering
\includegraphics[width=8cm]{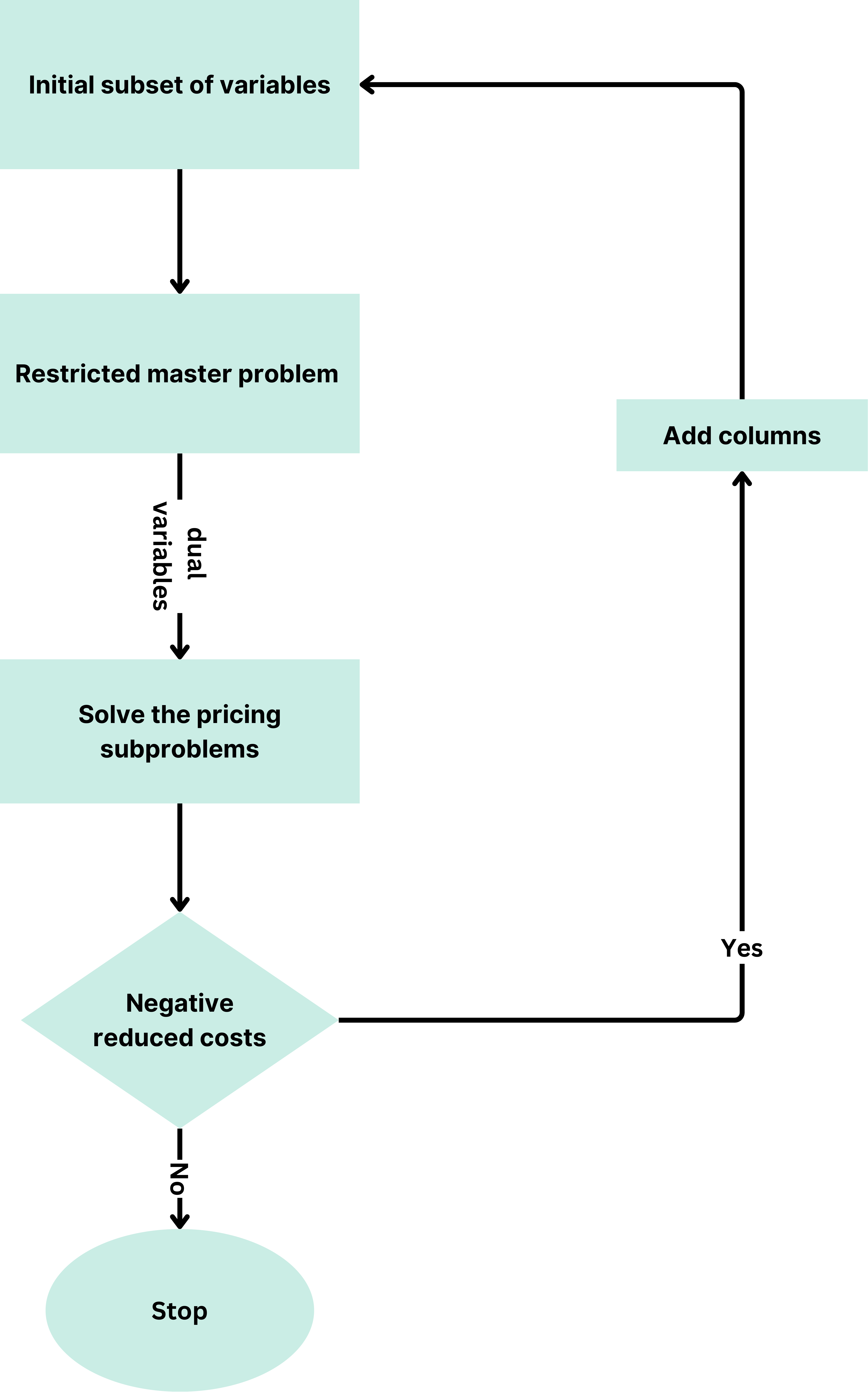}
\caption{Schematic for Dantzig Wolfe decomposition}
\label{fig:cg}
\end{figure}

The Dantzig-Wolfe algorithm reformulates the original problem by decomposing it into a master problem and a set of subproblems. In the Dantzig-Wolfe framework, the master problem initially contains only a subset of the variables, and new variables (columns) are iteratively generated by solving the subproblems. This process continues until an optimal solution is found which is implied by non-negative reduced costs, as shown in Figure \ref{fig:cg}.

 We now discuss Dantzig Wolfe decomposition in detail. Consider solving the problem \eqref{Decomp:FP}. The problem can be modified into \eqref{lag1dec}. The nonanticipativity constraints are complicating constraints. The structure problem now becomes like Figure \ref{fig:block:ld_cg}.  The problem can be further written as follows:
 
\begin{subequations}
      \begin{align}
\min \quad  &\sum_{s \in S}\left (\frac{c^T \bm{x}_s}{|\set{S}|}+z_s \right ) & \\
\text { s. t.} \quad & \bm{x}_1 = \bm{x}_s  \quad \forall s \in \set{S},s \neq 1 & \\
& \bm{x}_s,z_s \in X_s \quad \forall s \in \set{S} &
\end{align}
 \end{subequations}  

where \(X_s = \{\bm{x},z|A \bm{x}=b, \bm{x} \geq 0,T_{s} \bm{x}+W_{s} \bm{y}=h_{s},\bm{y} \geq 0,z = q^T_s \bm{y}  \}.\) 

Let \((\Tilde{\bm{x}}_{x,s}^j,\Tilde{\bm{x}}_{z,s}^j)\) be the extreme points and \((\Tilde{\bm{w}}_{x,s}^k,\Tilde{\bm{w}}_{z,s}^k)\) be extreme rays of \(X_s.\) The points in the polyhedrons \(X_s\) can be represented by a convex combination of the extreme points and a conic combination of the extreme rays of the polyhedron \citep{Bertsimas1997IntroductionOptimization}. Taking this into account, we can reformulate the LP in the following way:

\begin{subequations}
      \begin{align}
&\min \quad  \sum_{s \in S} \sum_{j \in J_s} \lambda_s^j\left (\frac{c^T \Tilde{\bm{x}}_{x,s}^j}{|\set{S}|}+\Tilde{x}_{z,s}^j \right ) + \sum_{s \in S} \sum_{k \in K_s} \theta_s^k\left (\frac{c^T \Tilde{\bm{w}}_{x,s}^k}{|\set{S}|}+\Tilde{w}_{z,s}^k \right ) \\
&\text { s. t.} \quad  \sum_{j \in J_1} \lambda_1^j \Tilde{\bm{x}}_{x,1}^j+\sum_{k \in K_1} \theta_1^k \Tilde{\bm{w}}_{x,1}^k = \sum_{j \in J_s} \lambda_s^j \Tilde{\bm{x}}_{x,s}^j+\sum_{k \in K_s} \theta_s^k \Tilde{\bm{w}}_{x,s}^k \quad \forall s \in \set{S},s \neq 1  \label{bpdual}\\
& \sum_{j \in J_s} \lambda_s^j = 1\quad \forall  s \in \set{S}   \label{bpdual1}\\ 
 &\lambda_s^j \geq 0,\quad \forall s \in \set{S},\;j\in J_{s}, \quad  \theta_s^k \geq 0\quad s \in \set{S},\;k\in K_{s} 
\end{align}
 \end{subequations} 
We call this reformulated problem the master problem. The master problem is the linear programming dual of the Lagrangian master problem \citep{Conforti2014ReformulationsRelaxations}.

Let \(\bm{\nu}_s\) be the set of dual variables for the constraint \eqref{bpdual} corresponding to \(s\) and \(r_s\) be the dual of the constraint \eqref{bpdual1} corresponding to \(s\).

To check if the solution obtained is optimal, we need to check the reduced costs of the variables and examine whether any one of them is negative. The reduced cost of \(\lambda_1^j\) is \(\frac{c^T \Tilde{\bm{x}}_{x,1}^j}{|\set{S}|}+\Tilde{x}_{z,1}^j-\sum_{s \in \set{S},s \neq 1}\bm{\nu}_s^T\Tilde{\bm{x}}_{x,1}^j-r_1\) and the reduced cost of \(\theta_1^k\) is  \(\frac{c^T \Tilde{\bm{w}}_{x,1}^k}{|\set{S}|}+\Tilde{w}_{z,1}^k-\sum_{s \in \set{S},s \neq 1}\bm{\nu}_s^T\Tilde{\bm{w}}_{x,1}^k.\) Similarly, for \(s \in S,s \neq 1,\) the reduced cost of \(\lambda_s^j\) is \(\frac{c^T \Tilde{\bm{x}}_{x,s}^j}{|\set{S}|}+\Tilde{x}_{z,s}^j+\bm{\nu}_s^T\Tilde{\bm{x}}_{x,s}^j-r_s\) and the reduced cost of \(\theta_s^k\) is  \(\frac{c^T \Tilde{\bm{w}}_{x,s}^k}{|\set{S}|}+\Tilde{w}_{z,s}^k+\bm{\nu}_s^T\Tilde{\bm{w}}_{x,s}^k.\)

While solving the complete master problem with all the extreme points and rays gives us the solution, it is very difficult to obtain all the extreme rays and points. Therefore, we start with a restricted master problem with a few extreme points and extreme rays and iteratively add them to the restricted master problem. 

 The goal is to identify if any reduced costs are negative. 

 The extreme points and extreme rays are columns that can be generated by solving the following pricing sub problems:

 \[
\begin{aligned}
\text{SP}_{1} : \min &\frac{c^T \bm{x}}{|\set{S}|}+z-\sum_{s \in \set{S},s \neq 1}\bm{\nu}_s^T\bm{x} \\
\text{s.t.} & \quad {\bm{x},z} \in X_1
\end{aligned}
\]

 \[
\begin{aligned}
\text{SP}_{s\in \set{S},s\neq 1} :  \min &\frac{c^T \bm{x}}{|\set{S}|}+z+\bm{\nu}_s^T\bm{x} \\
\text{s.t.} & \quad {\bm{x},z} \in X_s
\end{aligned}
\]

 If the pricing problem gives an objective of \(-\infty\), the solution is an extreme ray which can then be added to the master problem. If the pricing problem gives a finite objective which is less than the \(r_s\), the solution is an extreme point that can then be added to the master problem.  If the optimal objective from all the pricing subproblems \(SP_s\) are finite and greater than \(r_s,\) the reduced cost of all the variables \(\lambda_s^j,\theta_s^k\) are non-negative, implying that the solution obtained from the master problem is optimal. 

The steps for solving a linear programming problem using Dantzig Wolfe decomposition are shown in Algorithm \ref{alg:colgen}. 

\begin{algorithm}[H]
\caption{Dantzig Wolfe Decomposition Algorithm}
\label{alg:colgen}
\begin{algorithmic}[1]
\State \textbf{Initialize:} Start with a basic feasible solution to the master problem.
\While{True}
    \State Solve the master problem and obtain the dual variables for the constraints
    \State Solve the pricing subproblems $\text{SP}_s$.
    \If{all subproblems yield finite optimal objectives that are greater than correponding dual \(r_s\)}
        \State \textbf{Stop:} Current solution is optimal.
    \Else
        \State For subproblems returning an optimal objective of \(-\infty\) add the solution as an extreme ray to the master problem.
        \State For subproblems returning a finite optimal objective less than \(r_s\) add the solution as an extreme point to the master problem.
    \EndIf
    \State Repeat.
\EndWhile
\end{algorithmic}
\end{algorithm}

While the master problem gives upper bounds because of its restricted nature, summing up the pricing subproblems gives lower bounds due to its relaxed nature. Note that unlike Benders decomposition and Lagrangian decomposition, where the solution of the master problem converges from a smaller value to a larger one, the master problem in Dantzig-Wolfe decomposition moves in the opposite direction, from a larger value to a smaller one. This difference arises because in Dantzig-Wolfe, the master problem initially starts with a restricted set of variables and, with each iteration, adds more columns (subsets of variables), gradually expanding the feasible region. As more variables are introduced, the solution improves, leading to a lower objective value. In contrast, Benders and Lagrangian decompositions iteratively introduce new constraints or penalties in the form of cuts or relaxations, progressively tightening the feasible region of the master problem. As a result, their master problems start with a broader feasible region and converge upwards as they are constrained more heavily with each iteration.

The steps explained in a;gorithm \ref{alg:colgen} can be used directly to solve linear programming problems with complicating constraints. However, they cannot be directly used to solve  Mixed Integer Linear Programs (MILPs) with complicating constraints. In such cases, we can embed column generation methods like the Dantzig Wolfe decomposition algorithm, into the branch and bound framework. The branch and bound method is a method commonly used to solve MILPs \citep{Conforti2014GettingStarted}. The branch and bound method begins by solving the linear programming (LP) relaxation of the MILP, which ignores the integer constraints and provides an initial solution and bound on the objective value. If the solution to this relaxation includes fractional values for integer variables, the problem is split into child problems by branching on these fractional variables — creating new child problems with additional constraints that force the variable to be either rounded up or down. Each child problem is then solved, and bounds on the objective function are calculated. Problems that are infeasible or have bounds worse than the current best solution are pruned, eliminating them from further consideration. The algorithm selects the next new problem to explore based on strategies like depth-first or best-bound-first search and continues branching, bounding, and pruning until all child problems are either solved or discarded. 

To accelerate the branch and bound process, which involves solving multiple LPs, column generation methods can be integrated to solve the LPs more efficiently. This hybrid combination is known as branch and Price. Branch and price can also address some Mixed Integer Nonlinear Programming (MINLP) problems by constructing a Dantzig-Wolfe reformulation of the problem via a discretization approach \citep{Allman2021Branch-and-pricePrograms}. Application of branch and price to solve multi-time scale optimization models include the work by \cite{Gharaei2019AScheduling}.  In this study, \cite{Gharaei2019AScheduling} formulated a mixed-integer linear programming (MILP) model for an integrated production scheduling and distribution problem with routing decisions in a multi-site supply chain. The problem was addressed using a branch-and-price framework, with initial columns generated by a Bees algorithm (BA). The proposed framework was tested on several examples, demonstrating its effectiveness compared to commercial solvers and yielding tighter bounds compared to a stand-alone branch-and-price algorithm without the Bees algorithm.

\subsubsection{Limitations of decomposition algorithms}
An important limitation encountered while using decomposition algorithms is that these algorithms are not scalable to very large problems. While applying these algorithms to very large problems with say tens of millions of variables a large number of iterations are needed to converge and that could require a large amount of time. Another important aspect of these algorithms is that these decomposition algorithms are restricted to problems of a particular structure.

 \subsection{Metaheuristic algorithms} \label{metaheuristic}
 Metaheuristic algorithms are problem-independent methods often inspired by natural processes, designed to find approximate solutions to optimization problems by exploring a search space. They provide a framework that guides the process of searching for a solution toward optimal or near-optimal solutions by combining various heuristic methods and incorporating randomness and iterative improvements. This is done by balancing between exploring new regions of the search space (exploration) and intensifying the search around the best solutions found so far (exploitation).
 
 An example of a metaheuristic algorithm is the genetic algorithm (GA) \cite{Mitchell1996AnAlgorithms.} inspired by the process of natural selection. It involves the steps shown in 
Algorithm \ref{alg:genetic}.

\begin{algorithm}[H]
\caption{Genetic Algorithm}
\label{alg:genetic}
\begin{algorithmic}[1]
\State \textbf{Initialize:} Generate an initial population of solutions (chromosomes).
\While{stopping criterion not met}
    \State Evaluate the fitness/objective of each solution in the population.
    \State Select some of the solutions based on their fitness using a suitable heuristic \label{selection}
    \State Apply crossover and mutation operations to generate new solutions.
    \State Replace some or all of the old population with the newly generated solutions.
\EndWhile
\end{algorithmic}
\end{algorithm}
For line \ref{selection}, there have been a lot of heuristics described in literature like proportionate reproduction, ranking selection, and tournament selection \citep{Goldberg1991AAlgorithms}.

An application of the genetic algorithm to solve multi-time scale optimization problems includes the work by \cite{Lee2019SustainableRepresentation} who proposed a genetic algorithm to address integrated process planning and scheduling (IPPS) problems. IPPS problems are challenging combinatorial optimization problems characterized by numerous flexibilities/alternatives and constraints. These complexities make IPPS particularly difficult to solve. To effectively manage these challenges, \cite{Lee2019SustainableRepresentation} proposed a novel integrated chromosome representation that includes all the operations of an IPPS network and is capable of incorporating various flexibilities into a single string.  This new representation allows the use of a simple adaptation and application of the Genetic Algorithm. The proposed method was tested on a set of benchmark problems, and the results showed that it improved the makespan compared to other recently developed metaheuristics in significantly shorter computation times.

Other metaheuristic algorithms used to solve multi-time scale optimization models include simulated annealing  \Citep{Li2007AScheduling} inspired by the annealing process in metallurgy, ant colony optimization \Citep{Liu2016ApplicationScheduling} inspired by the actions of an ant colony, and variable neighborhood search \Citep{Leite2023SolvingAlgorithms} that involves systematically exploring the neighborhoods of solutions.  

Despite their versatility and effectiveness, metaheuristic methods face scalability issues, particularly when dealing with problems involving more than a thousand variables \citep{Hussain2019MetaheuristicSurvey}.
 
 \subsection{Matheuristic algorithms}\label{mathheuristic}
 Another method to solve multi-time scale optimization models involves algorithms that combine mathematical programming and heuristics which are called matheuristics \Citep{Ball2011HeuristicsProgramming}. An example of this is local branching \citep{Fischetti2003LocalBranching}, a technique used in mixed-integer programming (MIP) to improve solutions by iteratively exploring a sequence of solution neighborhoods in an MILP model. The steps involved in local branching are shown in Algorithm \ref{alg:local_branching}.

\begin{algorithm}[H]
\caption{Local Branching Algorithm}
\label{alg:local_branching}
\begin{algorithmic}[1]
\State Use a heuristic method to generate an initial feasible solution \( x^{\text{init}}\).
\State Set initial neighborhood size \( k \), and maximum iterations \( \text{max\_iter} \).
\State Initialize \( \text{iter} \gets 0 \) and \( x^{\text{best}} \gets x^{\text{init}} \).
\While{\( \text{iter} < \text{max\_iter} \)}
    \State Add local branching constraints to the MILP to limit the neighborhood to at most \( k \) changes from \( x^{\text{best}} \).
    \State Solve the modified MILP to obtain a solution \( x^{\text{new}} \).
    \If{an improved solution is found, i.e., \( f(x^{\text{new}}) < f(\text{best\_solution})\)}
        \State Update \( x^{\text{best}} \gets x^{\text{new}} \).
        \State Optionally reset or reduce \( k \) to intensify the search.
    \Else
        \State Increase \( k \) to explore a broader region.
    \EndIf
    \State Increment iteration counter: \( \text{iter} \gets \text{iter} + 1 \).
\EndWhile
\State Return \( \text{best\_solution} \).
\end{algorithmic}
\end{algorithm}

More complex and problem-specific matheuristics have been employed to solve multi-time scale optimization models. For instance, \cite{Ramanujam2023DistributedMicrogrid} proposed a mixed-integer linear programming (MILP) model that integrates the planning and scheduling of modular electrified plants, renewable-based generating units, and power lines within a microgrid. The decisions included single-time investment decisions on the location of the plants, renewable resources as well as power lines, monthly transportation and inventory decisions as well as hourly operating decisions. To solve the model, which included variables across three time scales, a K-means clustering-based aggregation-disaggregation matheuristic was introduced. This matheuristic involved two major steps: first, aggregating variables based on their locations into clusters using the K-means clustering algorithm on their coordinates, and second, disaggregating the investment decisions made into these clusters through a heuristic. To further obtain the power lines connecting locations in different clusters, a matching ILP was solved. Finally, to obtain the monthly and hourly decisions, the investment decisions were fixed in the full-space model. The matheuristic effectively reduced the problem's complexity by solving multiple small-sized MILPs instead of one large-sized MILP, making it particularly effective given the NP-hard nature of MILPs. The matheuristic was tested on a 20-location problem and was able to outperform commercial solvers.

 Other examples include the relax-and-fix heuristic proposed by \cite{Silva2023AProblem}, proposed to solve an integrated multiproduct, multiperiod, and multistage
capacitated lot sizing with a hybrid flow shop, as well as the spatial aggregation and decomposition algorithm proposed by \cite{Reinert2023ThisDecomposition} to solve large multi-time scale energy problems.

While matheuristics are efficient for problems with small to moderate size \Citep{Ramanujam2023DistributedMicrogrid}, they are not scalable to larger problems and are problem-specific. 
\subsection{Data-driven methods} \label{datadriven}
Data-driven methods leverage data to solve optimization models. There have been different data-driven methods used such as surrogate models and reinforcement-learning-based algorithms. The most common method used is the surrogate model which is a simpler model used to represent a complex model \Citep{Maravelias2009IntegrationOpportunities, Dias2020IntegrationModels,Badejo2022IntegratingAnalysis, Ye2015AUncertainty,Yang2023IntegratedResources, Chu2014IntegratedModeling,Beykal2022Data-drivenUncertainty}. An example of a surrogate-based method used to solve multi-time scale optimization models involves using surrogate models to represent the feasible region or carry out feasibility analysis \citep{Dias2019Data-drivenProblems}, particularly for the low-level variables. The surrogate model for the feasible low-level variables can then be included in a high-level model to obtain the high-level variables. The steps involved are shown in Algorithm \ref{alg:surrogate_model}.
\begin{algorithm}[H]
\caption{Surrogate Model Approach for High-Level Optimization}
\label{alg:surrogate_model}
\begin{algorithmic}[1]
\State \textbf{Step 1:} Formulate a high-level optimization model with unknown constraints for the low-level decision variables/ black box model.
\State \textbf{Step 2:} Obtain samples for the high-level decision variables and assess their output/feasibility in the low-level model.
\State \textbf{Step 3:} Determine algebraic approximations for the unknown constraints, through feasibility analysis and machine learning methods like regression and classification. These equations form a surrogate model for the detailed low-level model.
\State \textbf{Step 4:} Integrate the obtained constraints into the high-level model. Solve the high-level model to obtain good solutions.
\end{algorithmic}
\end{algorithm}

An application of this method can be found in the work by \cite{Dias2020IntegrationModels}, who utilized this approach to develop a framework for integrating planning, scheduling, and control. \cite{Dias2020IntegrationModels} initially focused on the integration of scheduling and control and treated the control problem as a black box (unknown constraints). A surrogate model was formulated to predict state and performance metrics and average flow based on control parameters. The integrated scheduling and control problem was then treated as a black box on the integration of planning and scheduling with a surrogate model formulated for ensuring the feasibility of aggregated production targets as well as predicting the energy consumption with respect to the integrated scheduling and control problem. To further manage the complexity and dimensionality of the problem, they introduced the concept of feature selection in the surrogate model construction, thereby streamlining the process and improving computational efficiency. The method was demonstrated by applying it to the optimization of an enterprise of air separation plants.
 
\begin{figure}[H]
\centering
\includegraphics[width=12cm]{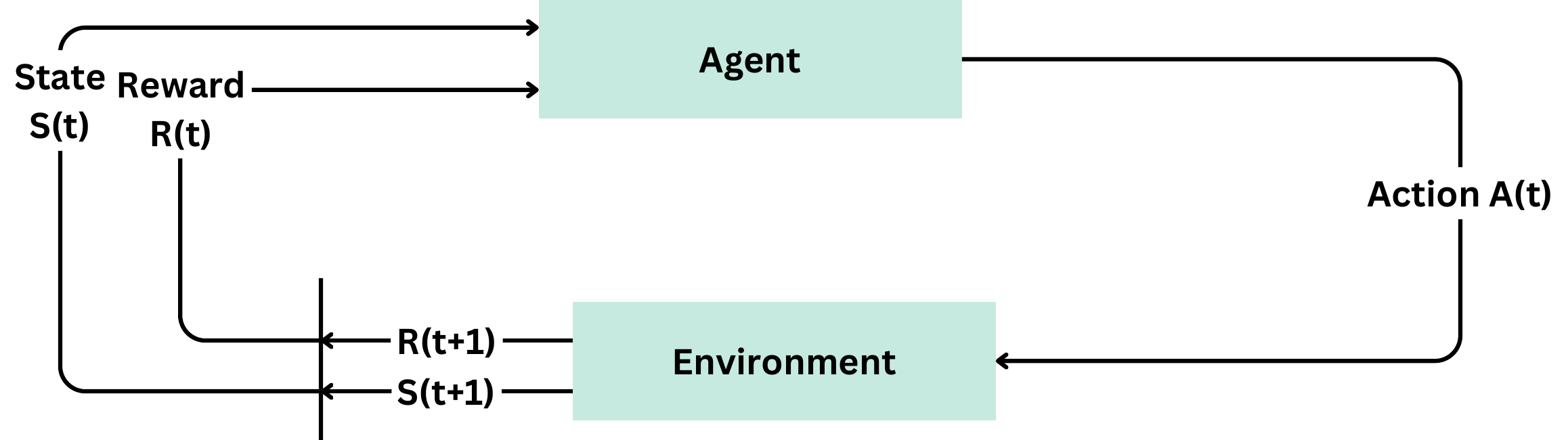}
\caption{Markov decision process}
\label{fig:mdp}
\end{figure}

Reinforcement learning-based approaches are another subset of data-driven methods used to solve multi-time scale optimization problems with multiple time periods \Citep{Shin2019Multi-timescaleProgramming,Ochoa2022Multi-agentMarkets}. Reinforcement learning (RL) algorithms learn optimal actions by interacting with an environment to maximize cumulative rewards. To implement RL algorithms for solving multi-time scale optimization models, the optimization model—or a component of it, such as the high-level version—is modeled as a Markov Decision Process (MDP) \citep{Sutton2018ReinforcementIntroduction}. In an MDP, at each time step, the process is in some state \(s\), and the decision maker may choose any action \(a\) that is available in state \(s\). The process responds at the next time step by moving into a new state \(s'\), and giving the decision maker a corresponding reward \(R_{a}(s,s').\) A representation of an MDP is shown in Figure \ref{fig:mdp}. The goal is to maximize the cumulative rewards. The MDP can then be solved using RL algorithms.

An example in the literature where this method was used is the work by \cite{Shin2019Multi-timescaleProgramming}, who developed a multi-time scale decision-making model that combines Markov decision process (MDP) and mathematical programming (MP) in a complementary way and introduced a computationally tractable solution algorithm based on reinforcement learning (RL) to solve the MP-embedded MDP problem. Additionally, \cite{Ochoa2022Multi-agentMarkets}  proposed a novel multi-agent deep reinforcement learning framework for efficient multi-time scale bidding for hybrid power plants. 

While these methods are effective for problems of moderate size, they struggle with scalability for larger problems due to the need for extensive parameter dimensions and a large number of samples.

 \subsection{PAMSO} \label{PAMSO}
\begin{figure}[H]
\centering
\includegraphics[width=15cm]{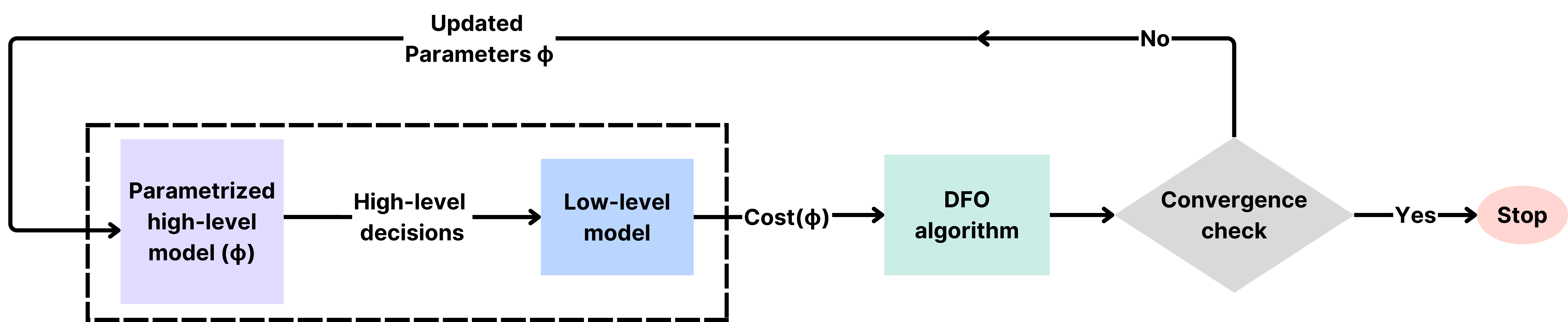}
\caption{Schematic for PAMSO}
\label{fig:PAMSO}
\end{figure}

 Parametric Autotuning Multi-time Scale Optimization (PAMSO) \citep{Ramanujam2024PAMSO:Algorithm} algorithm is a scalable and transferable algorithm used to solve multi-time scale optimization models and is inspired by the parametric cost function approximations (CFAs) \Citep{PerkinsIII2017StochasticApproximations} and the tuning of controllers using derivative-free optimizers  \Citep{Paulson2023ARepresentations, Coutinho2023BayesianControllers}. PAMSO involves tuning selected parameters in a low-fidelity optimization model to achieve a good performance in a high-fidelity optimization model as shown in Figure \ref{fig:PAMSO}.

PAMSO involves splitting the multi-time scale optimization model into a high-level model and low-level model. The high-level model provides high-level decisions which can then be fixed in the low-level model to obtain the low-level decisions as well as the objective of the multi-time scale optimization model. This is similar to the one-way communication method with \eqref{UP} being the high-level model and \eqref{LP} being the low-level model. 
  
We improve the performance of the high-level model by adding tunable parameters to the high-level model, to help better reflect the mismatch in variations between the high-level and low-level models. The parameterized high-level model is as follows:
\begin{subequations}
    \label{PUP}
    \begin{align}
        \min_{\bm{x},\bm{z}} \quad \Tilde{f}(\bm{x},\bm{\rho})+\Tilde{Q}(\bm{x},\bm{z},\bm{\Tilde{\theta}},\bm{\rho}) & \\
        \text{s.t.} \quad \Tilde{g}(\bm{x},\bm{\rho}) \leq 0 & \\
        \Tilde{H}(\bm{x}, \bm{z}, \bm{\Tilde{\theta}},\bm{\rho}) \leq 0 & 
    \end{align}
\end{subequations}

where \(\bm{\rho}\) are the set of tunable parameters we introduce. \(f,Q,g\) and \(H\) are modified into \(\Tilde{f},\Tilde{Q},\Tilde{g}\) and \(\Tilde{H}\),
respectively taking into account the parameters \(\bm{\rho}.\)

We now treat the system of high-level and low-level models as a black-box that takes in the tunable parameters as input and gives the corresponding objective of the multi-time scale optimization model as an output. This system is called the Multi-time scale black box function (MBBF). We can optimize the MBBF by tuning the parameters using Derivative Free (DFO) optimizers.  
The steps used in PAMSO are shown in Algorithm 
\ref{alg:PAMSO}.
\begin{algorithm}[H]
\caption{Optimization of Multi-Time Scale Model with Tunable Parameters}
\label{alg:PAMSO}
\begin{algorithmic}[1]
\State \textbf{Step 1:} Split the multi-time scale optimization model into a high-level model \eqref{UP} and a low-level model \eqref{LP}.
\State \textbf{Step 2:} Add tunable parameters to the high-level model.
\State \textbf{Step 3:} Treat the system comprising the parameterized high-level and low-level models as a black-box function. Optimize the black-box function by tuning the parameters using Derivative-Free Optimization (DFO) techniques.
\end{algorithmic}
\end{algorithm}

An important advantage of PAMSO is its ability to incorporate transfer learning by learning parameters from problems of smaller size to solve problems of large size. This can significantly reduce the time to solve large-sized problems. 

While PAMSO has been explored to solve case studies involving high-level decisions that are neither temporally interdependent nor variable over time, there is a need to check its effectiveness in cases where the high-level decisions are temporarily interdependent.

\section{Illustrative Example}
\label{sec:example}
In this section, we provide an example of a multi-time scale optimization problem and attempt to solve it with a few of the algorithms described in the previous section.
\subsection{Problem Statement}
An electric utility is installing two generators (indexed by $j=1,2$) with different fixed and operating
costs, in order to meet the demand within its service region for \textbf{3 representative days} \(s\). Each day is divided into three parts of 8 hours each, indexed by $i = 1,2,3$. These correspond to
parts of the day during which demand takes a base, medium, or peak value, respectively. The fixed cost per unit capacity of generator $j$ is amortized over its lifetime and amounts to $c_j$ per day. The operating cost of generator $j$ during the $i^{\text{th}}$ part of the day is $f_{i,j}$. If the demand during the $i^\text{th}$ part of the day cannot be served due to lack of capacity, additional capacity must be purchased at a cost of $g_{s}$ on the $s^{\text{th}}$ day per unit capacity purchased. 

The demand for the $i^{\text{th}}$ part of the day on the $s^{\text{th}}$ day  is $d_{i,s}.$ Furthermore, the availability of generator \(j\) on the $i^{\text{th}}$  part of the day for $s^{\text{th}}$  day  is $a_{s,i,j}.$ To optimize the cost, we need to optimize the installed capacity of the generator as well as the operating levels of the generators and the amount of power purchased to meet the unmet demand at different times of the day in different representative days.  For this, we formulate a multi-time scale model with decisions in 3-time scales: single-time decisions, and decisions for the different parts of different days (corresponding to 2 scales). The integrated model is formally defined in the next sub-section.

\subsection{Integrated model}
\subsubsection{Indices and Sets}
\begin{description}
    \item  [$j \in \set{J}$]  set of generators
    \item  [$i \in \set{I}$]  set of parts of the day for 8 hours each
    \item  [$s \in \set{S}$]  set of days
\end{description}
\subsubsection{Variables}
\begin{description}
\item  [$x_j$] installed capacity of generator $j$ (kW)
\item [$y_{s,i,j}$] operating levels of generator $j$ during the $i^\text{th}$ part of the day on the \(s^\text{th}\) day (kW)
\item [$\Tilde{y}_{s,i}$] power purchased in the $i^\text{th}$ part of the day on the $s^\text{th}$ day (kW) 
\end{description}
\subsubsection{Parameters}
\begin{description}
    \item[$a_{s,i,j}$] availability of generator $j$ on  $i^{\text{th}}$  part of the day for $s^{\text{th}}$  day 
    \item[$c_{j}$] fixed cost per unit capacity of generator $j$ per day amortized over its lifetime (\$/kW) 
    \item[$d_{s,i}$] demand of power during the $i^\text{th}$ part of the day on  the $s^{\text{th}}$ day (kW)
    \item[$f_{i,j}$] operating cost of generator $j$ during the $i^{\text{th}}$ part of the day (\$/kW)
    \item[$g_{s}$] cost of additional capacity purchased per power purchased on the $s^{\text{th}}$ day (\$/kW)
\end{description}
\subsubsection{Constraints}
We add a constraint to ensure that the capacity of each generator $j$ is  non-negative.

\begin{flalign}
   & x_{j} \geq 0 \quad \forall\: j \in \set{J}
   \label{minx}
\end{flalign}

The operating variables are non-negative. 

\begin{flalign}
   &  y_{s,i,j} \geq 0 \quad \forall\: i \in \set{I},s \in \set{S},j \in \set{J}
   \label{posy1}
\end{flalign}
\begin{flalign}
   & \Tilde{y}_{s,i}  \geq 0 \quad \forall\: i \in \set{I},s \in \set{S}
   \label{posy2}
\end{flalign}

The operating level of generators are constrained by $a_{s,i,j}x_{j}.$

\begin{flalign}
   & y_{s,i,j} \leq a_{s,i,j}x_{j} \quad \forall\: i \in \set{I},j \in \set{J},s \in \set{S}
   \label{avail}
\end{flalign}
The demand for power needs to be satisfied and is ensured by the following constraint:
\begin{flalign}
   & \sum_{j \in \set{J}}y_{s,i,j}+ \Tilde{y}_{s,i} \geq d_{s,i}\quad \forall\: i \in \set{I},k \in \set{K}
   \label{dem}
\end{flalign}
\subsubsection{Objective}

The objective of the model is to minimize the net cost consisting of the fixed cost of setting up the generators $\sum_{s \in \set{S}}\sum_{j \in \set{J}}c_{j}x_{j}$, as well as the total operating cost 
\\ $\sum_{s \in \set{S}}\left(\sum_{i\in \set{I}}\sum_{j \in \set{J}}\left(f_{i,j}y_{s,i,j} \right)+g_{s} \Tilde{y}_{s,i}\right)$. The objective is the sum of these components.
\subsubsection{Optimization model}
The entire model is given below:
\eqref{eq:fullmodel}
\begin{subequations} \label{eq:fullmodel}
    \begin{align}
        \min \quad & \phi =  \sum_{s \in \set{S}}\sum_{j \in \set{J}}c_{j}x_{j} + \sum_{s \in \set{S}}\left(\sum_{i\in \set{I}}\sum_{j \in \set{J}}\left(f_{i,j}y_{s,i,j} \right)+g_{s} \Tilde{y}_{s,i}\right) & \\
        \text{s.t. } \quad & x_{j} \geq 0 \quad \forall\: j \in \set{J} & \\
        & y_{s,i,j} \leq a_{s,i,j}x_{j} \quad \forall\: i \in \set{I},j \in \set{J},s \in \set{S} & \\
        & \sum_{j \in \set{J}}y_{s,i,j}+ \Tilde{y}_{s,i} \geq d_{s,i}\quad \forall\: i \in \set{I},s \in \set{S} & \\
         & y_{s,i,j} \geq 0 \quad \forall\: i \in \set{I},s \in \set{S},j \in \set{J} & \\
          &  \Tilde{y}_{s,i} \geq 0 \quad \forall\: i \in \set{I},s \in \set{S} &  
    \end{align}
\end{subequations}

\subsection{Solving the problem}
We now solve the problem using a few of the algorithms described in the previous section.  The algorithms are implemented in Python using pyomo models and the models are solved using Gurobi. The code for the implementation is given in \url{https://github.com/li-group/MultiScaleOpt-Tutorial.git}.
\subsubsection{Using full-space method}
We solve the multi-time scale optimization model and we obtain a cost of \$357,408.98 with \(x_1 = 2,515.15 \text{ kW},x_2 = 909.09 \text{ kW}.\)

\subsubsection{Using Benders Decomposition}
\begin{figure}[H]
\centering
\includegraphics[width=10cm]{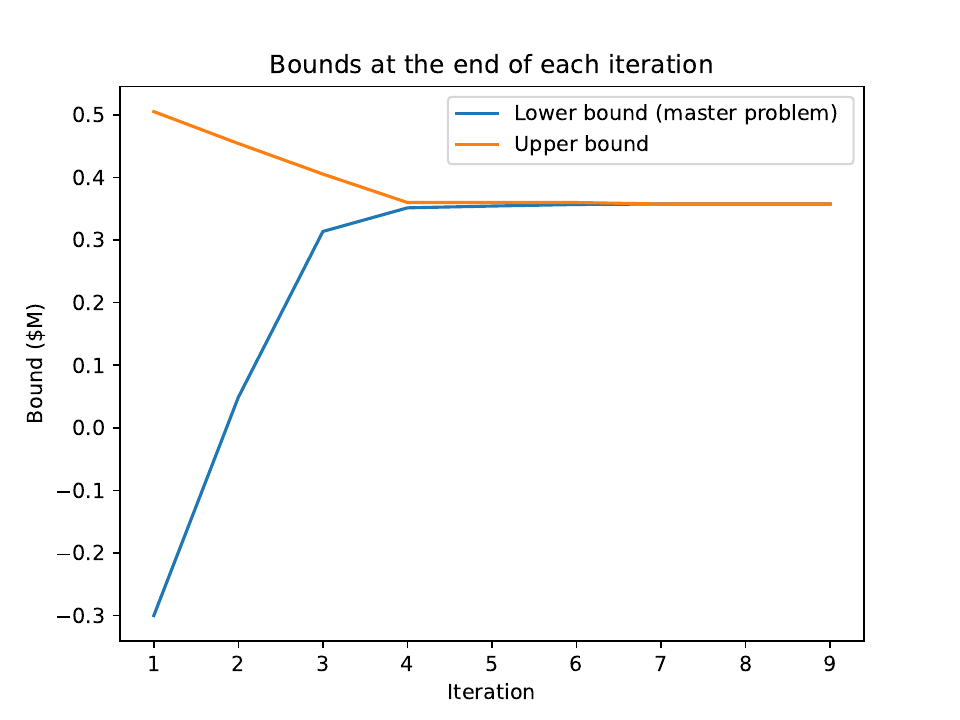}
\caption{Optimal objective of master problem obtained during Benders decomposition}
\label{fig:lbben}
\end{figure}

To implement the Benders decomposition, we formulate a master problem as follows:
\begin{subequations} \label{eq:hlmodel}
    \begin{align}
        \min \quad & \phi =  \sum_{s \in \set{S}}\sum_{j \in \set{J}}c_{j}x_{j} + \sum_{s \in \set{S}}z_s & \\
        \text{s.t. } \quad & x_{j} \geq 0 \quad \forall\: j \in \set{J} &  \\
        & z_s \geq -M \quad \forall\: s \in \set{S} \label{eq:lastconben} &
    \end{align}
\end{subequations}
where \(M\) is a large positive number. The constraint \eqref{eq:lastconben} is added to keep the master problem bounded.
We formulate subproblems as follows:
\begin{subequations} \label{eq:subfullmodel}
    \begin{align}
       SP_s =  \min \quad & z_s =  \sum_{j \in \set{J}}c_{j}x_{j} + \sum_{i\in \set{I}}\left(\sum_{j \in \set{J}}\left(f_{i,j}y_{i,j} \right)+g_{s} \Tilde{y}_{i}\right) & \\
        \text{s.t. } \quad & y_{i,j} \leq a_{s,i,j}x_{j}^* \quad \forall\: i \in \set{I},j \in \set{J} & \\
        & \sum_{j \in \set{J}}y_{i,j}+ \Tilde{y}_{i} \geq d_{s,i}\quad \forall\: i \in \set{I}& \\
         & y_{i,j} \geq 0 \quad \forall\: i \in \set{I},j \in \set{J} & \\
          &  \Tilde{y}_{i} \geq 0 \quad \forall\: i \in \set{I}&  
    \end{align}
\end{subequations}
where \(x^*\) is the solution obtained from the master problem.
We then use the benders decomposition algorithm as shown in section \ref{Benders} and get a solution of \(x_1 = 2515.15 \text{ kW}, x_2 = 909.09 \text{ kW}\) with a cost of \$357,408.97 which is the global minimum.
The algorithm runs for 9 iterations with the upper and lower bounds converging to the optimal solution as shown in figure \ref{fig:lbben}.
\subsubsection{Using Lagrangian decomposition}
To implement the Lagrangian decomposition, we decompose the problem into a set of subproblems as follows:
\begin{subequations}
      \begin{align}
 D_1({\nu}):= \min \quad &  \sum_{j \in \set{J}}c_{j}x_{j} + \sum_{i\in \set{I}}\left(\sum_{j \in \set{J}}\left(f_{i,j}y_{i,j} \right)+g_{1} \Tilde{y}_{i}\right)+\sum_{s \in \set{S},s \neq 1} \sum_{j}\nu_{s,j}x_j &  \\
\text { s. t.} \quad & x_j \geq 0  \quad \forall j & \\
& y_{i,j} \leq a_{1,i,j}x_{j} \quad \forall\: i \in \set{I},j \in \set{J},k \in \set{K} & \\
        & \sum_{j \in \set{J}}y_{i,j}+ \Tilde{y}_{i} \geq d_{1,i}\quad \forall\: i \in \set{I} & \\
         & y_{i,j} \geq 0 \quad \forall\: i \in \set{I},j \in \set{J} & \\
          &  \Tilde{y}_{i} \geq 0 \quad \forall\: i \in \set{I}, &  
\end{align}
 \end{subequations}

\begin{subequations}
      \begin{align}
 D_{s \in \set{S},s \neq 1}({\nu}):= \min \quad &  \sum_{j \in \set{J}}c_{j}x_{j} + \sum_{i\in \set{I}}\left(\sum_{j \in \set{J}}\left(f_{i,j}y_{i,j} \right)+g_{s} \Tilde{y}_{i}\right)-\sum_{j}\nu_{s,j}x_j &  \\
\text { s. t.} \quad & x_j \geq 0   \quad \forall j & \\
& y_{i,j} \leq a_{s,i,j}x_{j} \quad \forall\: i \in \set{I},j \in \set{J},k \in \set{K} & \\
        & \sum_{j \in \set{J}}y_{i,j}+ \Tilde{y}_{i} \geq d_{s,i}\quad \forall\: i \in \set{I} & \\
         & y_{i,j} \geq 0 \quad \forall\: i \in \set{I},j \in \set{J} & \\
          &  \Tilde{y}_{i} \geq 0 \quad \forall\: i \in \set{I} &  
\end{align}
 \end{subequations}

\begin{figure}[H]
\centering
\includegraphics[width=10cm]{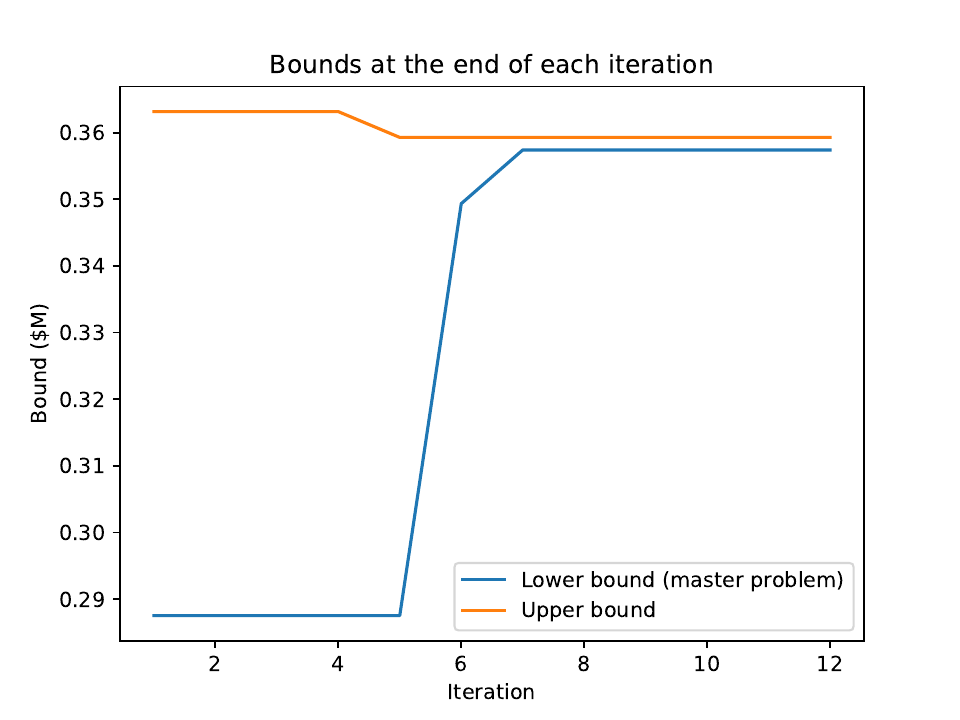}
\caption{Lower bound obtained during Lagrangian decomposition with cutting plane method}
\label{fig:lblagcp}
\end{figure}

We use the cutting plane method and the subgradient method to solve the problem using Lagrangian decomposition. The heuristic we use to update the upper bound and obtain feasible solutions is to fix all the investment decisions obtained from solving the subproblems at various iterations.
The cutting plane method gives a lower bound of \(\$357,408.98\) which is the global minimum and an upper bound of \(\$ 359,290.31\) which is around 0.5\% from the global optimum with \(x_1 =1876.21 \text{ kW},x_2 = 1619.03 \text{ kW}.\) The decomposition algorithm runs for 12 iterations till we obtain the same Lagrangian multipliers. The upper and lower bounds at the end of each iteration are shown in Figure \ref{fig:lblagcp}.

\begin{figure}[H]
\centering
\includegraphics[width=10cm]{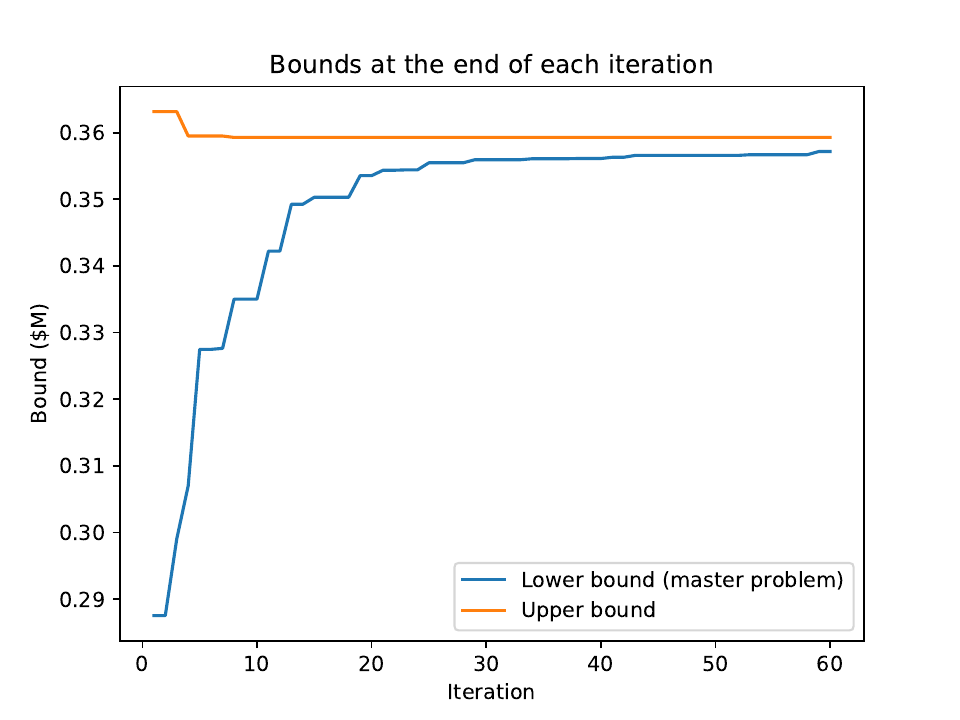}
\caption{Lower bound obtained during Lagrangian decomposition with subgradient method}
\label{fig:lblagsg}
\end{figure}

We run the subgradient method for 60 iterations. We get a lower bound of \(\$357,169.21\) which is around approximately the global optimum and an upper bound of \(\$ 359,290.31\) which is around 0.5\% from the global optimum with \(x_1 =1876.21 \text{ kW},x_2 = 1619.03 \text{ kW}.\) The upper and lower bounds obtained at the end of each iteration are shown in Figure \ref{fig:lblagsg}.

\subsubsection{Using Dantzig Wolfe decomposition}
To implement  Dantzig Wolfe decomposition to the problem, we formulate a master problem as follows:

\begin{subequations}
      \begin{align}
&\min \quad  \sum_{s \in S} \sum_{k \in K_{\lambda,s}} \lambda_s^k\left (\sum_{j \in \set{J}}c_{j}\Tilde{x}_{x,s,j}^k+\Tilde{x}_{z,s}^k \right ) + \sum_{s \in S} \sum_{k \in K_{\theta,s}} \theta_s^k\left (\sum_{j \in \set{J}}c_{j}\Tilde{w}_{x,s,j}^k+\Tilde{w}_{z,s}^k \right ) \\
&\text { s. t.} \quad  \sum_{k \in K_{\lambda,1}} \lambda_1^k \Tilde{x}_{x,1,j}^k+\sum_{k \in K_{\theta,1}} \theta_1^k \Tilde{w}_{x,1,j} =\sum_{k \in K_{\lambda,a}}\lambda_s^k \Tilde{x}_{x,s,j}+\sum_{k \in K_{\theta,s}} \theta_s^k \Tilde{w}_{x,1,j} \quad \forall s \in \set{S},s \neq 1  \label{bpdual_imp}\\
& \sum_{k \in K_{\lambda,s}} \lambda_s^k = 1\quad \forall  s \in \set{S}   \label{bpdual1_imp}\\ 
 &\lambda_s^k \geq 0,\quad \forall k \in S_{\lambda,s},\;s \in \set{S}, \quad  \theta_s^k \geq 0\quad k \in S_{\theta,s},\;s \in \set{S}
\end{align}
 \end{subequations} 

We then formulate the following subproblems as follows:

 \[
\begin{aligned}
\text{SP}_{1} : \min &\sum_{j \in \set{J}}c_{j}x_{j}+z-\sum_{j \in \set{J}}\sum_{s \in \set{S},s \neq 1}\nu_{s,j}x_j \\
\text{s.t.} \quad &  z =  \sum_{j \in \set{J}}c_{j}x_{j} + \sum_{i\in \set{I}}\left(\sum_{j \in \set{J}}\left(f_{i,j}y_{i,j} \right)+g_{1} \Tilde{y}_{i}\right) & \\
         & y_{i,j} \leq a_{1,i,j}x_{j} \quad \forall\: i \in \set{I},j \in \set{J} & \\
        & \sum_{j \in \set{J}}y_{i,j}+ \Tilde{y}_{i} \geq d_{1,i}\quad \forall\: i \in \set{I}& \\
         & y_{i,j} \geq 0 \quad \forall\: i \in \set{I},j \in \set{J} & \\
          &  \Tilde{y}_{i} \geq 0 \quad \forall\: i \in \set{I} & \\
          & x_j \geq 0 \quad \forall\; j\in \set{J} & \\
           & x_j \leq M \quad \forall\; j\in \set{J}
\end{aligned}
\]

 \[
\begin{aligned}
\text{SP}_{s\in \set{S},s\neq 1} :  \min &\sum_{j \in \set{J}}c_{j}x_{j}+z+\sum_{j \in \set{J}}\nu_{s,j}x \\
\text{s.t.} \quad &  z =  \sum_{j \in \set{J}}c_{j}x_{j} + \sum_{i\in \set{I}}\left(\sum_{j \in \set{J}}\left(f_{i,j}y_{i,j} \right)+g_{s} \Tilde{y}_{i}\right) & \\
         & y_{i,j} \leq a_{s,i,j}x_{j} \quad \forall\: i \in \set{I},j \in \set{J} & \\
        & \sum_{j \in \set{J}}y_{i,j}+ \Tilde{y}_{i} \geq d_{s,i}\quad \forall\: i \in \set{I}& \\
         & y_{i,j} \geq 0 \quad \forall\: i \in \set{I},j \in \set{J} & \\
          &  \Tilde{y}_{i} \geq 0 \quad \forall\: i \in \set{I} & \\
          & x_j \geq 0 \quad \forall\; j\in \set{J} & \\
          & x_j \leq M \quad \forall\; j\in \set{J}
\end{aligned}
\]
where \(\nu_{s,j}\) is the dual for the constraint \eqref{bpdual_imp} corresponding to \(s\) for the \(j^\text{th}\) generator.

\begin{figure}[H]
\centering
\includegraphics[width=10cm]{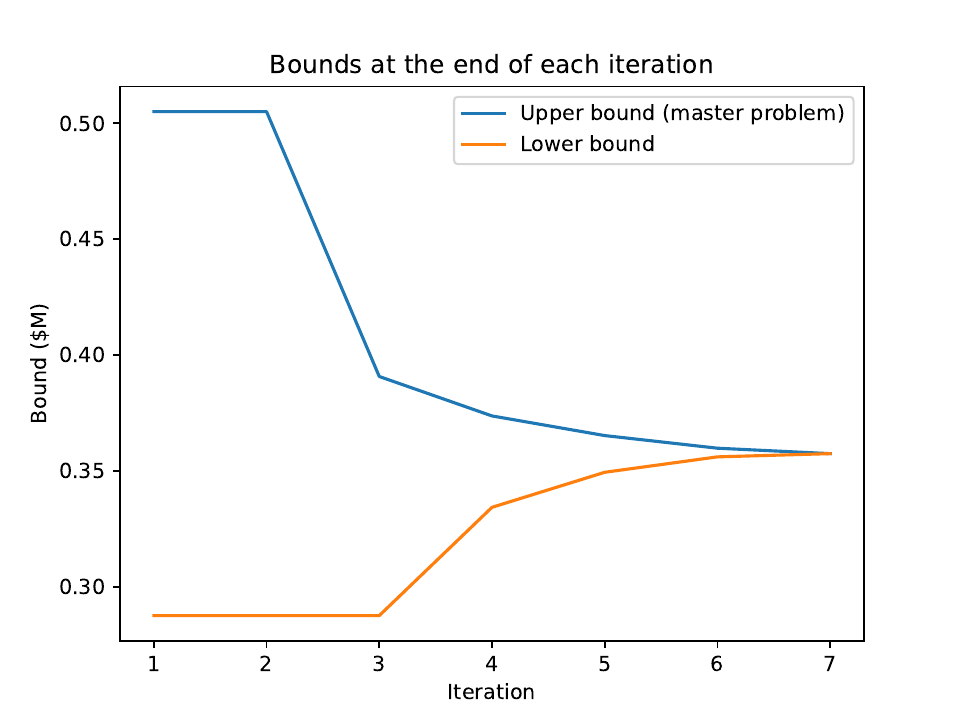}
\caption{Optimal objective of master problem obtained during Dantzig Wolfe decomposition}
\label{fig:ubdw}
\end{figure}

We add constraints upper bounding the investment decisions to keep the subproblem bounded.
We apply the Dantzig Wolfe decomposition algorithm with the above formulations. The algorithm runs for 7 iterations before getting all non-negative reduced costs and we get a solution of \(\bm{x}_1 = 2,515.15, x_2 = 909.09\) with a cost of \$357,408.97 which is the global minimum. The algorithm runs for 7 iterations with the upper and lower bounds converging to the optimal solution as shown in figure \ref{fig:ubdw}.

\subsubsection{Using PAMSO}
 To implement PAMSO, we split the optimization model into a high-level model taking decisions regarding the capacity of the 2 generators, and a low-level model deciding the operating levels of generators as well as the power purchased. 

 The high-level model is adapted from the full space model with the operating decisions for the entire duration aggregated. The high-level model takes decisions in only one time scale, i.e., decisions valid for the entire time duration as a whole. With a slight abuse in notation, the operating variables of the aggregated model are summarized below:

\begin{table}
\caption{Variables in high-level model \label{varhlmodel}}{%
\begin{tabular}{c c c c}
\toprule
Variable in  & Physical  & Variable in  & Physical  \\
full-space & meaning & high-level model & meaning \\
\midrule
        $y_{s,i,j}$ & Operating levels & $y_{j}$ & Net operating \\
         & of generator $j$  &  & levels of generator  \\
         &  during the $i^\text{th}$ & & generator $j$ \\
         &  part of the  & &   \\
         &  $s^\text{th}$ day &  & \\
        \hline
        $\Tilde{y}_{s,i}$ & Power purchased & $\Tilde{y}$ & Net power \\
        & during the $i^\text{th}$ &  &  purchased   \\
        & part of  &  &   \\
         & the $s^\text{th}$ day   &  &  \\
\botrule
\end{tabular}}{}
\end{table}

From the high-level model, we obtain the installed capacity of the generators \(x_j.\)
The optimization model is given below:
\begin{subequations} \label{eq:hlmodelex}
    \begin{align}
        \min \quad & \phi =  \sum_{s \in \set{S}}\sum_{j \in \set{J}}c_{j}x_{j} + \sum_{s \in \set{S}}\left(\sum_{j \in \set{J}}\sum_{i\in \set{I}}\left(f_{i,j}\right)y_{j} +g_{s} \Tilde{y}\right) & \\
        \text{s.t. } \quad & x_{j} \geq 0 \quad \forall\: j \in \set{J} & \\
        & y_{j} \leq \sum_{s \in \set{S}}\sum_{i \in \set{I}}a_{s,i,j}x_{j} \quad \forall\: ,j \in \set{J} & \\
        & \sum_{j \in \set{J}}y_{j}+ \Tilde{y} \geq  \sum_{s \in \set{S}}\sum_{i \in \set{I}}d_{s,i} &\\
        & y_{j} \geq 0 \quad \forall\:j \in \set{J} & \\
          &  \Tilde{y} \geq 0 &  
    \end{align}
\end{subequations}

The low-level model is the full-space model with the installed capacity \(x_{j}\) of the generators fixed from the high-level model. The low-level model takes decisions for the parts of the day on the 3 days.

\textbf{Parameterized high-level model:} Based on the data, we can either underestimate or overestimate the effect of the availability in the aggregated model. Therefore we add a prefactor to the availability of each of the generators to compensate for the mismatch.   Furthermore, to better capture the differences in the operation of the generators, we add a third parameter which signifies the minimum capacity of each generator i.e., the capacity of each generator is greater than this parameter. This can help in cases when an optimum configuration consists of having both generators installed with a reasonable capacity. Without this parameter, there is a high chance of installing only one generator for cases where there are some intricate differences between the generators. Therefore, the parameters we add to the high-level model are
\begin{enumerate}
    \item \(\rho_1\) which represents the prefactor to the availability of generator 1,
    \item \(\rho_2\) which represents the prefactor to the availability of generator 2,
    \item \(\rho_3\) which represents the minimum capacity  installed for each of the generators.
\end{enumerate}

\begin{subequations} \label{eq:paramhlmodel}
    \begin{align}
        \min \quad & \phi =  \sum_{s \in \set{S}}\sum_{j \in \set{J}}c_{j}x_{j} + \sum_{s \in \set{S}}\left(\sum_{j \in \set{J}}\sum_{i\in \set{I}}\left(f_{i,j}\right)y_{j} +g_{s} \Tilde{y}\right) & \\
        \text{s.t. } \quad & x_{j} \geq 0 \quad \forall\: j \in \set{J} & \\
        & y_j \leq \rho_j\sum_{s \in \set{S}}\sum_{i \in \set{I}}a_{i,k,j}x_j \quad \forall\: j \in \set{J}  & \\
        & \sum_{j \in \set{J}}y_{j}+ \Tilde{y} \geq  \sum_{k \in \set{K}}\sum_{i \in \set{I}}d_{i,k} & \\
        & x_{j} \geq \rho_3 \quad \forall\: j \in \set{J} & \\
        & y_{j} \geq 0 \quad \forall\:j \in \set{J} & \\
          &  \Tilde{y} \geq 0 &  
    \end{align}
\end{subequations}

We then form the MBBF using the parameterized high-level model and the low-level model and tune the MBBF using DFO solvers.

\begin{figure}[H]
\centering
\includegraphics[width=10cm]{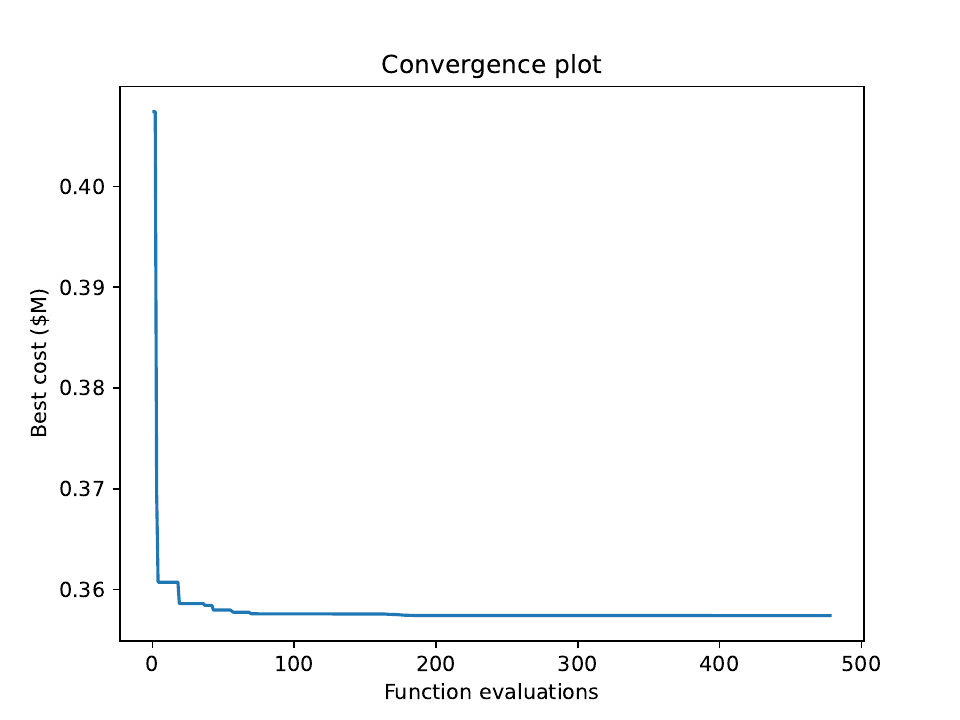}
\caption{Best cost obtained in different function evaluations of PAMSO}
\label{fig:conPAMSO}
\end{figure}

For our implementation of PAMSO, we use the Bayesian Adaptive Direct Search (BADS) \citep{Acerbi2017PracticalSearch} algorithm as a DFO solver via the PyBADS software \citep{Singh2024PyBADS:Python} to solve the MBBF. The range for the parameters is \(\rho_1 \in [0,1.5],\rho_2 \in [0,1.5], \rho_3 \in [0,2100].\) After applying the BADS algorithm we get the optimum parameters as \(\rho_1 = 0.609,\rho_2 = 0.279, \rho_3 = 909.092\) which gives a solution of \(x_1 = 2515.15 \text{ kW},x_2 = 909.09 \text{ kW} \) and a cost of \$357,408.98 which is the global minimum. The algorithm takes 478 function evaluations to reach the solution. A convergence plot, showing the evolution of the best cost with function evaluations is shown in Figure \ref{fig:conPAMSO}.

\subsection{VMM}
We finally calculate the VMM of the problem with the high-level model described in Equations \eqref{eq:hlmodelex}. The optimal objective from solving the multi-time scale model (MM) turns out to be \$357,408.98 . The Multi-scale Performance of the Single-scale Solution (MPSS) is \(\$407,520.75.\) Thus the VMM is \$50,111.77.

\section{Conclusion} \label{sec:conclusion}
 In this tutorial, we have provided an overview of multi-time scale optimization models and the algorithms that have been used in the literature to solve these models. Multi-timescale optimization models involve decision variables in multiple time scales like long-time/high-level planning decisions and short-time/low-level scheduling and control
decisions. The applications of these models are widespread,  including in economics and environmental management. The \textit{Value of the Multi-scale Model} (VMM) is a metric introduced to quantify the benefit of using a multi-time scale optimization model as opposed to using separate independent models for the different decision layers. 

Multi-time scale optimization models can be of different sizes and different structures. To solve such models, different methods have been used in the literature. We have provided the different methods, their applications, and their limitations. We have also applied some of the algorithms to solve an example optimization problem at \url{https://github.com/li-group/MultiScaleOpt-Tutorial.git}.  
 
 \bibliography{references.bib}

\begin{thebibliography}{53}
\providecommand{\natexlab}[1]{#1}
\providecommand{\url}[1]{\texttt{#1}}
\expandafter\ifx\csname urlstyle\endcsname\relax
  \providecommand{\doi}[1]{doi: #1}\else
  \providecommand{\doi}{doi: \begingroup \urlstyle{rm}\Url}\fi

\bibitem[Acerbi and Ma(2017)]{Acerbi2017PracticalSearch}
Luigi Acerbi and Wei~Ji Ma.
\newblock {Practical Bayesian Optimization for Model Fitting with Bayesian Adaptive Direct Search}.
\newblock In I~Guyon, U~Von Luxburg, S~Bengio, H~Wallach, R~Fergus, S~Vishwanathan, and R~Garnett, editors, \emph{Advances in Neural Information Processing Systems}, volume~30. Curran Associates, Inc., 2017.
\newblock URL \url{https://proceedings.neurips.cc/paper_files/paper/2017/file/df0aab058ce179e4f7ab135ed4e641a9-Paper.pdf}.

\bibitem[Allen et~al.(2023)Allen, Baratsas, Kakodkar, Avraamidou, Demirhan, Heuberger-Austin, Klokkenburg, and Pistikopoulos]{Allen2023ASystems}
Richard~Cory Allen, Stefanos~G Baratsas, Rahul Kakodkar, Styliani Avraamidou, Cosar~Doga Demirhan, Clara~F Heuberger-Austin, Mark Klokkenburg, and Efstratios~N Pistikopoulos.
\newblock {A multi-period integrated planning and scheduling approach for developing energy systems}.
\newblock \emph{Optimal Control Applications and Methods}, 44\penalty0 (2):\penalty0 355--372, 3 2023.
\newblock ISSN 0143-2087.
\newblock \doi{https://doi.org/10.1002/oca.2866}.
\newblock URL \url{https://doi.org/10.1002/oca.2866}.

\bibitem[Allman and Zhang(2021)]{Allman2021Branch-and-pricePrograms}
Andrew Allman and Qi~Zhang.
\newblock {Branch-and-price for a class of nonconvex mixed-integer nonlinear programs}.
\newblock \emph{Journal of Global Optimization}, 81\penalty0 (4):\penalty0 861--880, 2021.
\newblock ISSN 1573-2916.
\newblock \doi{10.1007/s10898-021-01027-w}.
\newblock URL \url{https://doi.org/10.1007/s10898-021-01027-w}.

\bibitem[Badejo and Ierapetritou(2022)]{Badejo2022IntegratingAnalysis}
Oluwadare Badejo and Marianthi Ierapetritou.
\newblock {Integrating tactical planning, operational planning and scheduling using data-driven feasibility analysis}.
\newblock \emph{Computers {\&} Chemical Engineering}, 161:\penalty0 107759, 2022.
\newblock ISSN 0098-1354.
\newblock \doi{https://doi.org/10.1016/j.compchemeng.2022.107759}.
\newblock URL \url{https://www.sciencedirect.com/science/article/pii/S0098135422001004}.

\bibitem[Ball(2011)]{Ball2011HeuristicsProgramming}
Michael~O Ball.
\newblock {Heuristics based on mathematical programming}.
\newblock \emph{Surveys in Operations Research and Management Science}, 16\penalty0 (1):\penalty0 21--38, 2011.
\newblock ISSN 1876-7354.
\newblock \doi{https://doi.org/10.1016/j.sorms.2010.07.001}.
\newblock URL \url{https://www.sciencedirect.com/science/article/pii/S1876735410000036}.

\bibitem[Barzanji et~al.(2020)Barzanji, Naderi, and Begen]{Barzanji2020DecompositionProblem}
Ramin Barzanji, Bahman Naderi, and Mehmet~A Begen.
\newblock {Decomposition algorithms for the integrated process planning and scheduling problem}.
\newblock \emph{Omega}, 93:\penalty0 102025, 2020.
\newblock ISSN 0305-0483.
\newblock \doi{https://doi.org/10.1016/j.omega.2019.01.003}.
\newblock URL \url{https://www.sciencedirect.com/science/article/pii/S0305048318306698}.

\bibitem[Bertsekas(1997)]{Bertsekas1997NonlinearProgramming}
Dimitri~P Bertsekas.
\newblock {Nonlinear programming}.
\newblock \emph{Journal of the Operational Research Society}, 48\penalty0 (3):\penalty0 334, 1997.

\bibitem[Bertsimas and Tsitsiklis(1997)]{Bertsimas1997IntroductionOptimization}
Dimitris Bertsimas and John Tsitsiklis.
\newblock \emph{{Introduction to Linear Optimization}}.
\newblock Athena Scientific, 1st edition, 1997.
\newblock ISBN 1886529191.

\bibitem[Beykal et~al.(2022)Beykal, Avraamidou, and Pistikopoulos]{Beykal2022Data-drivenUncertainty}
Burcu Beykal, Styliani Avraamidou, and Efstratios~N Pistikopoulos.
\newblock {Data-driven optimization of mixed-integer bi-level multi-follower integrated planning and scheduling problems under demand uncertainty}.
\newblock \emph{Computers {\&} Chemical Engineering}, 156:\penalty0 107551, 2022.
\newblock ISSN 0098-1354.
\newblock \doi{https://doi.org/10.1016/j.compchemeng.2021.107551}.
\newblock URL \url{https://www.sciencedirect.com/science/article/pii/S009813542100329X}.

\bibitem[Biondi et~al.(2017)Biondi, Sand, and Harjunkoski]{Biondi2017OptimizationApproach}
Matteo Biondi, Guido Sand, and Iiro Harjunkoski.
\newblock {Optimization of multipurpose process plant operations: A multi-time-scale maintenance and production scheduling approach}.
\newblock \emph{Computers {\&} Chemical Engineering}, 99:\penalty0 325--339, 2017.
\newblock ISSN 0098-1354.
\newblock \doi{https://doi.org/10.1016/j.compchemeng.2017.01.007}.
\newblock URL \url{https://www.sciencedirect.com/science/article/pii/S009813541730008X}.

\bibitem[Birge and Louveaux(2011)]{Birge2011IntroductionProgramming}
John~R. Birge and François Louveaux.
\newblock \emph{{Introduction to Stochastic Programming}}.
\newblock Springer New York, New York, NY, 2011.
\newblock ISBN 978-1-4614-0236-7.
\newblock \doi{10.1007/978-1-4614-0237-4}.

\bibitem[Chu and You(2013)]{Chu2013IntegratedApproach}
Yunfei Chu and Fengqi You.
\newblock {Integrated Scheduling and Dynamic Optimization of Complex Batch Processes with General Network Structure Using a Generalized Benders Decomposition Approach}.
\newblock \emph{Industrial {\&} Engineering Chemistry Research}, 52\penalty0 (23):\penalty0 7867--7885, 6 2013.
\newblock ISSN 0888-5885.
\newblock \doi{10.1021/ie400475s}.
\newblock URL \url{https://doi.org/10.1021/ie400475s}.

\bibitem[Chu and You(2014)]{Chu2014IntegratedModeling}
Yunfei Chu and Fengqi You.
\newblock {Integrated Planning, Scheduling, and Dynamic Optimization for Batch Processes: MINLP Model Formulation and Efficient Solution Methods via Surrogate Modeling}.
\newblock \emph{Industrial {\&} Engineering Chemistry Research}, 53\penalty0 (34):\penalty0 13391--13411, 8 2014.
\newblock ISSN 0888-5885.
\newblock \doi{10.1021/ie501986d}.
\newblock URL \url{https://doi.org/10.1021/ie501986d}.

\bibitem[Conforti et~al.(2014{\natexlab{a}})Conforti, Cornu{\'{e}}jols, and Zambelli]{Conforti2014GettingStarted}
Michele Conforti, Gérard Cornu{\'{e}}jols, and Giacomo Zambelli.
\newblock {Getting Started}.
\newblock In Michele Conforti, Gérard Cornu{\'{e}}jols, and Giacomo Zambelli, editors, \emph{Integer Programming}, pages 1--44. Springer International Publishing, Cham, 2014{\natexlab{a}}.
\newblock ISBN 978-3-319-11008-0.
\newblock \doi{10.1007/978-3-319-11008-0{\_}1}.
\newblock URL \url{https://doi.org/10.1007/978-3-319-11008-0_1}.

\bibitem[Conforti et~al.(2014{\natexlab{b}})Conforti, Cornu{\'{e}}jols, and Zambelli]{Conforti2014ReformulationsRelaxations}
Michele Conforti, Gérard Cornu{\'{e}}jols, and Giacomo Zambelli.
\newblock {Reformulations and Relaxations}.
\newblock In Michele Conforti, Gérard Cornu{\'{e}}jols, and Giacomo Zambelli, editors, \emph{Integer Programming}, pages 321--350. Springer International Publishing, Cham, 2014{\natexlab{b}}.
\newblock ISBN 978-3-319-11008-0.
\newblock \doi{10.1007/978-3-319-11008-0{\_}8}.
\newblock URL \url{https://doi.org/10.1007/978-3-319-11008-0_8}.

\bibitem[Coutinho et~al.(2023)Coutinho, Santos, and Reis]{Coutinho2023BayesianControllers}
João P~L Coutinho, Lino~O Santos, and Marco~S Reis.
\newblock {Bayesian Optimization for automatic tuning of digital multi-loop PID controllers}.
\newblock \emph{Computers {\&} Chemical Engineering}, 173:\penalty0 108211, 2023.
\newblock ISSN 0098-1354.
\newblock \doi{https://doi.org/10.1016/j.compchemeng.2023.108211}.
\newblock URL \url{https://www.sciencedirect.com/science/article/pii/S0098135423000807}.

\bibitem[Dias and Ierapetritou(2019)]{Dias2019Data-drivenProblems}
Lisia~S. Dias and Marianthi~G. Ierapetritou.
\newblock {Data-driven feasibility analysis for the integration of planning and scheduling problems}.
\newblock \emph{Optimization and Engineering}, 20\penalty0 (4):\penalty0 1029--1066, 12 2019.
\newblock ISSN 1389-4420.
\newblock \doi{10.1007/s11081-019-09459-w}.

\bibitem[Dias and Ierapetritou(2020)]{Dias2020IntegrationModels}
Lisia~S Dias and Marianthi~G Ierapetritou.
\newblock {Integration of planning, scheduling and control problems using data-driven feasibility analysis and surrogate models}.
\newblock \emph{Computers {\&} Chemical Engineering}, 134:\penalty0 106714, 2020.
\newblock ISSN 0098-1354.
\newblock \doi{https://doi.org/10.1016/j.compchemeng.2019.106714}.
\newblock URL \url{https://www.sciencedirect.com/science/article/pii/S0098135419306982}.

\bibitem[Dowling et~al.(2017)Dowling, Kumar, and Zavala]{Dowling2017AParticipation}
Alexander~W Dowling, Ranjeet Kumar, and Victor~M Zavala.
\newblock {A multi-scale optimization framework for electricity market participation}.
\newblock \emph{Applied Energy}, 190:\penalty0 147--164, 2017.
\newblock ISSN 0306-2619.
\newblock \doi{https://doi.org/10.1016/j.apenergy.2016.12.081}.
\newblock URL \url{https://www.sciencedirect.com/science/article/pii/S0306261916318487}.

\bibitem[Erdirik-Dogan and Grossmann(2008)]{Erdirik-Dogan2008SimultaneousLines}
Muge Erdirik-Dogan and Ignacio~E Grossmann.
\newblock {Simultaneous planning and scheduling of single-stage multi-product continuous plants with parallel lines}.
\newblock \emph{Computers {\&} Chemical Engineering}, 32\penalty0 (11):\penalty0 2664--2683, 2008.
\newblock ISSN 0098-1354.
\newblock \doi{https://doi.org/10.1016/j.compchemeng.2007.07.010}.
\newblock URL \url{https://www.sciencedirect.com/science/article/pii/S0098135407001998}.

\bibitem[Fischetti and Lodi(2003)]{Fischetti2003LocalBranching}
Matteo Fischetti and Andrea Lodi.
\newblock {Local branching}.
\newblock \emph{Mathematical Programming}, 98\penalty0 (1):\penalty0 23--47, 2003.
\newblock ISSN 1436-4646.
\newblock \doi{10.1007/s10107-003-0395-5}.
\newblock URL \url{https://doi.org/10.1007/s10107-003-0395-5}.

\bibitem[Gharaei and Jolai(2019)]{Gharaei2019AScheduling}
Ali Gharaei and Fariborz Jolai.
\newblock {A branch and price approach to the two-agent integrated production and distribution scheduling}.
\newblock \emph{Computers {\&} Industrial Engineering}, 136:\penalty0 504--515, 2019.
\newblock ISSN 0360-8352.
\newblock \doi{https://doi.org/10.1016/j.cie.2019.07.050}.
\newblock URL \url{https://www.sciencedirect.com/science/article/pii/S0360835219304474}.

\bibitem[Goldberg and Deb(1991)]{Goldberg1991AAlgorithms}
David~E Goldberg and Kalyanmoy Deb.
\newblock {A Comparative Analysis of Selection Schemes Used in Genetic Algorithms}.
\newblock In GREGORY J~E RAWLINS, editor, \emph{Foundations of Genetic Algorithms}, volume~1, pages 69--93. Elsevier, 1991.
\newblock ISBN 1081-6593.
\newblock \doi{https://doi.org/10.1016/B978-0-08-050684-5.50008-2}.
\newblock URL \url{https://www.sciencedirect.com/science/article/pii/B9780080506845500082}.

\bibitem[Guignard and Kim(1987)]{Guignard1987LagrangeanBounds}
Monique Guignard and Siwhan Kim.
\newblock {Lagrangean decomposition: A model yielding stronger lagrangean bounds}.
\newblock \emph{Mathematical Programming}, 39\penalty0 (2):\penalty0 215--228, 1987.
\newblock ISSN 1436-4646.
\newblock \doi{10.1007/BF02592954}.
\newblock URL \url{https://doi.org/10.1007/BF02592954}.

\bibitem[Hooker and Ottosson(2003)]{Hooker2003Logic-basedDecomposition}
John~N Hooker and Greger Ottosson.
\newblock {Logic-based Benders decomposition}.
\newblock \emph{Mathematical Programming}, 96\penalty0 (1):\penalty0 33--60, 2003.

\bibitem[Hussain et~al.(2019)Hussain, Mohd~Salleh, Cheng, and Shi]{Hussain2019MetaheuristicSurvey}
Kashif Hussain, Mohd~Najib Mohd~Salleh, Shi Cheng, and Yuhui Shi.
\newblock {Metaheuristic research: a comprehensive survey}.
\newblock \emph{Artificial Intelligence Review}, 52\penalty0 (4):\penalty0 2191--2233, 2019.
\newblock ISSN 1573-7462.
\newblock \doi{10.1007/s10462-017-9605-z}.
\newblock URL \url{https://doi.org/10.1007/s10462-017-9605-z}.

\bibitem[Kim et~al.(2024)Kim, Choi, Park, Adams, Heo, and Lee]{Kim2024Multi-periodUncertainty}
Sunwoo Kim, Yechan Choi, Joungho Park, Derrick Adams, Seongmin Heo, and Jay~H Lee.
\newblock {Multi-period, multi-timescale stochastic optimization model for simultaneous capacity investment and energy management decisions for hybrid Micro-Grids with green hydrogen production under uncertainty}.
\newblock \emph{Renewable and Sustainable Energy Reviews}, 190:\penalty0 114049, 2024.
\newblock ISSN 1364-0321.
\newblock \doi{https://doi.org/10.1016/j.rser.2023.114049}.
\newblock URL \url{https://www.sciencedirect.com/science/article/pii/S1364032123009073}.

\bibitem[Kopanos and Puigjaner(2009)]{Kopanos2009Multi-SiteIndustries}
Georgios~M Kopanos and Luis Puigjaner.
\newblock {Multi-Site Scheduling/Batching and Production Planning for Batch Process Industries}.
\newblock In Rita~Maria de~Brito~Alves, Caludio Augusto~Oller do~Nascimento, and Evaristo~Chalbaud Biscaia, editors, \emph{Computer Aided Chemical Engineering}, volume~27, pages 2109--2114. Elsevier, 2009.
\newblock ISBN 1570-7946.
\newblock \doi{https://doi.org/10.1016/S1570-7946(09)70742-4}.
\newblock URL \url{https://www.sciencedirect.com/science/article/pii/S1570794609707424}.

\bibitem[Lee and Ha(2019)]{Lee2019SustainableRepresentation}
Hyun~Cheol Lee and Chunghun Ha.
\newblock {Sustainable Integrated Process Planning and Scheduling Optimization Using a Genetic Algorithm with an Integrated Chromosome Representation}.
\newblock \emph{Sustainability}, 11\penalty0 (2), 2019.
\newblock ISSN 2071-1050.
\newblock \doi{10.3390/su11020502}.
\newblock URL \url{https://www.mdpi.com/2071-1050/11/2/502}.

\bibitem[Leite et~al.(2023)Leite, Pinto, and Alves]{Leite2023SolvingAlgorithms}
Mário Manuel~Silva Leite, Telmo Miguel~Pires Pinto, and Cláudio Manuel~Martins Alves.
\newblock {Solving the integrated planning and scheduling problem using variable neighborhood search based algorithms}.
\newblock \emph{Expert Systems with Applications}, 228:\penalty0 120191, 10 2023.
\newblock ISSN 0957-4174.
\newblock \doi{10.1016/J.ESWA.2023.120191}.

\bibitem[Li et~al.(2022)Li, Conejo, Liu, Omell, Siirola, and Grossmann]{Li2022Mixed-integerSystems}
Can Li, Antonio~J Conejo, Peng Liu, Benjamin~P Omell, John~D Siirola, and Ignacio~E Grossmann.
\newblock {Mixed-integer linear programming models and algorithms for generation and transmission expansion planning of power systems}.
\newblock \emph{European Journal of Operational Research}, 297\penalty0 (3):\penalty0 1071--1082, 2022.
\newblock ISSN 0377-2217.
\newblock \doi{https://doi.org/10.1016/j.ejor.2021.06.024}.
\newblock URL \url{https://www.sciencedirect.com/science/article/pii/S0377221721005397}.

\bibitem[Li and McMahon(2007)]{Li2007AScheduling}
W~D Li and C~A McMahon.
\newblock {A simulated annealing-based optimization approach for integrated process planning and scheduling}.
\newblock \emph{International Journal of Computer Integrated Manufacturing}, 20\penalty0 (1):\penalty0 80--95, 1 2007.
\newblock ISSN 0951-192X.
\newblock \doi{10.1080/09511920600667366}.
\newblock URL \url{https://doi.org/10.1080/09511920600667366}.

\bibitem[Liu et~al.(2016)Liu, Ni, and Qiu]{Liu2016ApplicationScheduling}
Xiaojun Liu, Zhonghua Ni, and Xiaoli Qiu.
\newblock {Application of ant colony optimization algorithm in integrated process planning and scheduling}.
\newblock \emph{The International Journal of Advanced Manufacturing Technology}, 84\penalty0 (1):\penalty0 393--404, 2016.
\newblock ISSN 1433-3015.
\newblock \doi{10.1007/s00170-015-8145-4}.
\newblock URL \url{https://doi.org/10.1007/s00170-015-8145-4}.

\bibitem[Maravelias and Sung(2009)]{Maravelias2009IntegrationOpportunities}
Christos~T Maravelias and Charles Sung.
\newblock {Integration of production planning and scheduling: Overview, challenges and opportunities}.
\newblock \emph{Computers {\&} Chemical Engineering}, 33\penalty0 (12):\penalty0 1919--1930, 2009.
\newblock ISSN 0098-1354.
\newblock \doi{https://doi.org/10.1016/j.compchemeng.2009.06.007}.
\newblock URL \url{https://www.sciencedirect.com/science/article/pii/S0098135409001501}.

\bibitem[Mitchell(1996)]{Mitchell1996AnAlgorithms.}
Melanie Mitchell.
\newblock \emph{{An introduction to genetic algorithms.}}
\newblock Complex adaptive systems. The MIT Press, Cambridge, MA, US, 1996.
\newblock ISBN 0-262-13316-4 (Hardcover).

\bibitem[Munawar et~al.(2005)Munawar, Kapadi, Patwardhan, Madhavan, Pragathieswaran, Lingathurai, and Gudi]{Munawar2005IntegrationManufacturing}
S~A Munawar, Mangesh~D Kapadi, S~C Patwardhan, K~P Madhavan, S~Pragathieswaran, P~Lingathurai, and Ravindra~D Gudi.
\newblock {Integration of planning and scheduling in multi-site plants: Application to paper manufacturing}.
\newblock In Luis Puigjaner and Antonio Espu{\~{n}}a, editors, \emph{Computer Aided Chemical Engineering}, volume~20, pages 1621--1626. Elsevier, 2005.
\newblock ISBN 1570-7946.
\newblock \doi{https://doi.org/10.1016/S1570-7946(05)80112-9}.
\newblock URL \url{https://www.sciencedirect.com/science/article/pii/S1570794605801129}.

\bibitem[Ochoa et~al.(2022)Ochoa, Gil, Angulo, and Valle]{Ochoa2022Multi-agentMarkets}
Tomás Ochoa, Esteban Gil, Alejandro Angulo, and Carlos Valle.
\newblock {Multi-agent deep reinforcement learning for efficient multi-timescale bidding of a hybrid power plant in day-ahead and real-time markets}.
\newblock \emph{Applied Energy}, 317:\penalty0 119067, 2022.
\newblock ISSN 0306-2619.
\newblock \doi{https://doi.org/10.1016/j.apenergy.2022.119067}.
\newblock URL \url{https://www.sciencedirect.com/science/article/pii/S0306261922004603}.

\bibitem[Paulson et~al.(2023)Paulson, Sorourifar, and Mesbah]{Paulson2023ARepresentations}
Joel~A. Paulson, Farshud Sorourifar, and Ali Mesbah.
\newblock {A Tutorial on Derivative-Free Policy Learning Methods for Interpretable Controller Representations}.
\newblock In \emph{2023 American Control Conference (ACC)}, pages 1295--1306. IEEE, 5 2023.
\newblock ISBN 979-8-3503-2806-6.
\newblock \doi{10.23919/ACC55779.2023.10156412}.

\bibitem[Peng et~al.(2019)Peng, Zhang, Feng, Rong, and Su]{Peng2019AUncertainty}
Zedong Peng, Yi~Zhang, Yiping Feng, Gang Rong, and Hongye Su.
\newblock {A Progressive Hedging-Based Solution Approach for Integrated Planning and Scheduling Problems under Demand Uncertainty}.
\newblock \emph{Industrial {\&} Engineering Chemistry Research}, 58\penalty0 (32):\penalty0 14880--14896, 8 2019.
\newblock ISSN 0888-5885.
\newblock \doi{10.1021/acs.iecr.9b02620}.
\newblock URL \url{https://doi.org/10.1021/acs.iecr.9b02620}.

\bibitem[Perkins~III and Powell(2017)]{PerkinsIII2017StochasticApproximations}
Raymond~T Perkins~III and Warren~B Powell.
\newblock {Stochastic optimization with parametric cost function approximations}.
\newblock \emph{arXiv preprint arXiv:1703.04644}, 2017.

\bibitem[Ramanujam and Li(2024)]{Ramanujam2024PAMSO:Algorithm}
Asha Ramanujam and Can Li.
\newblock {PAMSO: Parametric Autotuning Multi-time Scale Optimization Algorithm}.
\newblock 7 2024.

\bibitem[Ramanujam et~al.(2023)Ramanujam, Constante-Flores, and Li]{Ramanujam2023DistributedMicrogrid}
Asha Ramanujam, Gonzalo~E Constante-Flores, and Can Li.
\newblock {Distributed manufacturing for electrified chemical processes in a microgrid}.
\newblock \emph{AIChE Journal}, n/a\penalty0 (n/a):\penalty0 e18265, 10 2023.
\newblock ISSN 0001-1541.
\newblock \doi{https://doi.org/10.1002/aic.18265}.
\newblock URL \url{https://doi.org/10.1002/aic.18265}.

\bibitem[Reinert et~al.(2023)Reinert, Nilges, Baumg{\"{a}}rtner, and Bardow]{Reinert2023ThisDecomposition}
Christiane Reinert, Benedikt Nilges, Nils Baumg{\"{a}}rtner, and André Bardow.
\newblock {This is SpArta: Rigorous Optimization of Regionally Resolved Energy Systems by Spatial Aggregation and Decomposition}.
\newblock \emph{arXiv e-prints}, page arXiv:2302.05222, 2 2023.
\newblock \doi{10.48550/arXiv.2302.05222}.

\bibitem[Rockafellar and Wets(1991)]{Rockafellar1991ScenariosUncertainty}
R~T Rockafellar and Roger J.-B. Wets.
\newblock {Scenarios and Policy Aggregation in Optimization Under Uncertainty}.
\newblock \emph{Mathematics of Operations Research}, 16\penalty0 (1):\penalty0 119--147, 2 1991.
\newblock ISSN 0364-765X.
\newblock \doi{10.1287/moor.16.1.119}.
\newblock URL \url{https://doi.org/10.1287/moor.16.1.119}.

\bibitem[Shah and Ierapetritou(2012)]{Shah2012IntegratedIndustry}
Nikisha~K Shah and Marianthi~G Ierapetritou.
\newblock {Integrated production planning and scheduling optimization of multisite, multiproduct process industry}.
\newblock \emph{Computers {\&} Chemical Engineering}, 37:\penalty0 214--226, 2012.
\newblock ISSN 0098-1354.
\newblock \doi{https://doi.org/10.1016/j.compchemeng.2011.08.007}.
\newblock URL \url{https://www.sciencedirect.com/science/article/pii/S0098135411002626}.

\bibitem[Shin and Lee(2019)]{Shin2019Multi-timescaleProgramming}
Joohyun Shin and Jay~H Lee.
\newblock {Multi-timescale, multi-period decision-making model development by combining reinforcement learning and mathematical programming}.
\newblock \emph{Computers {\&} Chemical Engineering}, 121:\penalty0 556--573, 2019.
\newblock ISSN 0098-1354.
\newblock \doi{https://doi.org/10.1016/j.compchemeng.2018.11.020}.
\newblock URL \url{https://www.sciencedirect.com/science/article/pii/S0098135418308913}.

\bibitem[Silva and Mateus(2023)]{Silva2023AProblem}
D~M Silva and G~R Mateus.
\newblock {A Mixed-Integer Programming Formulation and Heuristics for an Integrated Production Planning and Scheduling Problem}.
\newblock In Luca Di~Gaspero, Paola Festa, Amir Nakib, and Mario Pavone, editors, \emph{Metaheuristics}, pages 290--305, Cham, 2023. Springer International Publishing.
\newblock ISBN 978-3-031-26504-4.

\bibitem[Singh and Acerbi(2024)]{Singh2024PyBADS:Python}
Gurjeet~Sangra Singh and Luigi Acerbi.
\newblock {PyBADS: Fast and robust black-box optimization in Python}.
\newblock \emph{Journal of Open Source Software}, 9\penalty0 (94):\penalty0 5694, 2 2024.
\newblock ISSN 2475-9066.
\newblock \doi{10.21105/joss.05694}.

\bibitem[Subramanian et~al.(2013)Subramanian, Rawlings, Maravelias, Flores-Cerrillo, and Megan]{Subramanian2013IntegrationManagement}
Kaushik Subramanian, James~B Rawlings, Christos~T Maravelias, Jesus Flores-Cerrillo, and Lawrence Megan.
\newblock {Integration of control theory and scheduling methods for supply chain management}.
\newblock \emph{Computers {\&} Chemical Engineering}, 51:\penalty0 4--20, 2013.
\newblock ISSN 0098-1354.
\newblock \doi{https://doi.org/10.1016/j.compchemeng.2012.06.012}.
\newblock URL \url{https://www.sciencedirect.com/science/article/pii/S0098135412001883}.

\bibitem[Sutton and Barto(2018)]{Sutton2018ReinforcementIntroduction}
Richard~S Sutton and Andrew~G Barto.
\newblock \emph{{Reinforcement learning: An introduction}}.
\newblock MIT press, 2018.

\bibitem[Tabrizi(2018)]{Tabrizi2018IntegratedImpacts}
Babak~H Tabrizi.
\newblock {Integrated planning of project scheduling and material procurement considering the environmental impacts}.
\newblock \emph{Computers {\&} Industrial Engineering}, 120:\penalty0 103--115, 2018.
\newblock ISSN 0360-8352.
\newblock \doi{https://doi.org/10.1016/j.cie.2018.04.031}.
\newblock URL \url{https://www.sciencedirect.com/science/article/pii/S0360835218301670}.

\bibitem[Yang et~al.(2023)Yang, Gao, and You]{Yang2023IntegratedResources}
Shiyu Yang, H~Oliver Gao, and Fengqi You.
\newblock {Integrated optimization in operations control and systems design for carbon emission reduction in building electrification with distributed energy resources}.
\newblock \emph{Advances in Applied Energy}, page 100144, 2023.
\newblock ISSN 2666-7924.
\newblock \doi{https://doi.org/10.1016/j.adapen.2023.100144}.
\newblock URL \url{https://www.sciencedirect.com/science/article/pii/S2666792423000239}.

\bibitem[Ye and You(2015)]{Ye2015AUncertainty}
Wenhe Ye and Fengqi You.
\newblock {A simulation-based optimization method for solving the integrated supply chain network design and inventory control problem under uncertainty}.
\newblock In \emph{Chemical Engineering Transactions}, volume~45, pages 499--504. Italian Association of Chemical Engineering - AIDIC, 10 2015.
\newblock \doi{10.3303/CET1545084}.
\newblock URL \url{http://www.scopus.com/inward/record.url?scp=84946116663&partnerID=8YFLogxK http://www.scopus.com/inward/citedby.url?scp=84946116663&partnerID=8YFLogxK}.

\end{thebibliography}

\end{document}